\documentclass[12pt]{article}
\usepackage{amsmath}
\usepackage{graphicx,psfrag,epsf}
\usepackage{enumerate}
\usepackage{natbib}
\usepackage{url} % not crucial - just used below for the URL

\newcommand{\blind}{0}

\usepackage{mathrsfs}
\usepackage{amsmath}
\usepackage{amssymb}
\usepackage{amsthm}
\usepackage{pifont}
\usepackage{color}
\usepackage{txfonts}
\usepackage{cases}
\usepackage{dsfont}
\usepackage{titlesec}
\usepackage{fancyhdr}
\usepackage{enumitem}
\usepackage{subfig}
\usepackage{bm,graphicx,amsfonts,makeidx}
\usepackage{comment}

\newtheorem{thm}{Theorem}[section]

\begin{document}

\def\spacingset#1{\renewcommand{\baselinestretch}%
{#1}\small\normalsize} \spacingset{1}

%%%%%%%%%%%%%%%%%%%%%%%%%%%%%%%%%%%%%%%%%%%%%%%%%%%%%%%%%%%%%%%%%%%%%%%%%%%%%%

\if0\blind
{
  \title{\bf Sparse Wavelet Estimation in Quantile Regression with Multiple Functional Predictors}
  \author{Dengdeng Yu\footnotemark[1], Li Zhang\thanks{These authors contributed equally.}, Ivan Mizera, Bei Jiang, and Linglong Kong\thanks{Corresponding author.
\textit{E-mail address}:  lkong@ualberta.ca (L. Kong)} \thanks{
    This
work has been supported by the \textit{Natural Sciences and Engineering Research
Council of Canada} and \textit{Canadian Statistical Sciences Institute}. } \\
%    and \\
%    Author 2 \\
   Department of Mathematical and Statistical
Sciences, \\
University of Alberta, \\
Edmonton, Alberta, Canada
    }
  \maketitle\thispagestyle{empty}
} \fi

\if1\blind
{
  \bigskip
  \bigskip
  \bigskip
  \begin{center}
    {\LARGE\bf Title}
\end{center}
  \medskip
} \fi

\bigskip
\begin{abstract} %200 or fewer words
In this manuscript, we study quantile regression in partial functional linear model where response is scalar and predictors include both scalars and multiple functions.
Wavelet basis are adopted to better approximate functional slopes while effectively detect local features.
The sparse group lasso penalty is imposed to select important functional predictors while capture shared information among them.
%This approach can be easily extended to the setting where the functional predictors are located in a two- or three- dimensional space.
The estimation problem can be reformulated into a standard second-order cone program and then solved by an interior point method. %implemented in a powerful R package: Rmosek.
We also give a novel algorithm by using alternating direction method of multipliers (ADMM) which was recently employed by many researchers in solving penalized quantile regression problems. The asymptotic properties such as the convergence rate and prediction error bound have been established. Simulations and a real data from ADHD-200 fMRI data are investigated to show the superiority of our proposed method.
%By utilizing wavelet bases,
%derive the asymptotic properties of consistency
\end{abstract}

\noindent%
{\it Keywords:}  Functional data analysis; Sparse group lasso; ADMM; Convergence rate; Prediction error bound;  ADHD
%3 to 6 keywords, that do not appear in the title
\vfill

%\footnotetext{ \textit{Preprint submitted to Journal of Computational and Graphical Statistics, June 5, 2017}}

\newpage
\spacingset{1.45} % DON'T change the spacing!

\section{Introduction}
Functional data analysis (FDA) is about the analysis of information on curves, images, functions, or more general objects.
It has become a major branch of nonparametric statistics and is a fast evolving area as more data has arisen where the primary object of observation can be viewed as a function (\citealp{ramsay2006functional}; \citealp{wang2015review}; \citealp{morris2015functional}).
A standard functional linear model with scalar response and functional covariate is
\begin{equation}\label{flm}
y=\alpha+\int_{0}^{1} x(t)\beta(t)dt +\varepsilon,
\end{equation}
where the coefficient $\beta(t)$ is a function, and $\varepsilon$ is a random error.
To estimate the functional coefficient $\beta(t)$, we can use functional basis to approximate it.
There are three major choices of functional basis:
 general basis such as B-spline basis and wavelet basis (\citealp{cardot2003spline}; \citealp{zhao2012wavelet}),
 functional principal component basis  (\citealp{Cardot1999}; \citealp{cai2006prediction}; \citealp{mulleryao2008}; %\citealp{delaigle2009}, %  \citealp{hall2007methodology},
 \citealp{kong2016partially}), and partial least square basis (\citealp{delaigle2012methodology}).
 Recently in imaging analysis, \cite{zhao2012wavelet}, \cite{wang2014regularized} and \cite{zhao2015wavelet} successfully adopted wavelet basis with regularizations to estimate
 the functional slope where the functional covariates are image features located in 1D, 2D and 3D domains respectively.

%, such as the roughness penalty approach with restricted wavelet basis function (\citealp{zhao2012wavelet}, \citealp{wang2014regularized}, \citealp{zhao2015wavelet}), can be implemented to produce an estimator of $\beta(t)$ that is meaningful in interpretation and useful in prediction.

The functional linear model (\ref{flm}) can be extended to a partial functional linear model with multiple functional covariates
\begin{equation}\label{mflm}
y=\alpha+ \int_0^1 \bm x^T(t)\bm \beta(t)dt+\bm u^{T}\bm{\gamma}+\varepsilon,
\end{equation}
where covariates $\bm u$ are scalars and $\bm \gamma$ are the coefficients.
The functional coefficients $\bm \beta(t)$ can be estimated by using regularization techniques.
In particular, penalized principal component basis has been an especially popular choice (\citealp{gertheiss2013variable}; \citealp{lian2013shrinkage}).
Recently, \cite{kong2016partially} successfully applied such technique to model (\ref{mflm}) in the setting of ultrahigh-dimensional scalar predictors.
%analyze the spectrometric data and Japanese weather data while \citet{gertheiss2013variable} used the  to investigate the magnetic resonance imaging (MRI) data to %select the important MRI indices for multiple sclerosis  disease prediction.

In recent years, quantile regression, which was introduced by the seminal work of \citet{koenker1978regression},
has been well developed and recognized in functional linear regression,
with many mainly focusing on the functional linear quantile regression model:
\begin{equation}\label{quantile:model0}
Q_{\tau}(y|x(t))=\alpha_{\tau}+\int_{0}^{1} x(t) \beta_{\tau}(t)dt,
\end{equation}
where $Q_{\tau}(y|x(t))$ is the $\tau$-th conditional quantile of response $y$ given a functional covariate $ x(t)$ for a fixed quantile level $\tau \in (0,1)$.
As an alternative to least squares regression, % in model (\ref{flm}),
the quantile regression method % in model (\ref{quantile:model0})
%is sometimes considered to be a more flexible and robust method with certain desirable properties. For instance, quantile regression,
is more efficient and robust when the responses are non-normal, errors are heavy tailed or outliers are present.
It is also capable of dealing with the heteroscedasticity issues and providing a more complete picture of the response
(\citealp{koenker2005quantile}).
%By estimating conditional quantile functions at several different level, % compared with mean regression,
%one may better account for the location, scale and shape of the response distribution.
%The monograph by \citet{koenker2005quantile} provides an excellent summary of the history and recent development in quantile regression.
%As in mean regression for simple functional linear model,
To estimate the functional coefficient $\beta_\tau(t)$,
functional basis can as well be used to approximate it; for instance,
%it is convenient to restrict it in a functional space with a finite basis,
general basis like B-spline basis (\citealp{cardot2005quantile}; \citealp{sun2005semiparametric}), functional principle component basis (\citealp{kato2012estimation}; \citealp{lu2014functional}; \citealp{tang2014partial}) and partial quantile basis (\citealp{yu2016partial}).
%As pointed out by \citealp{yao2017regularized}, quantile regression is particularly a motivated choice for health scores and economic indices which tend to be non-normal with severe skewness and outliers, and considered to be a robust estimation technique.

In this article, we extend model (\ref{quantile:model0})
to a partial functional linear quantile regression model with multiple functional covariates
%a more flexible setting with multiple functional covariates together with additional scalar vector such that
\begin{equation}\label{quantile:model}
Q_{\tau}(y|~\bm u,\bm x(t))=\alpha_{\tau}+\int_{0}^{1}\bm x^T(t)\bm \beta_{\tau}(t)dt+\bm u^T\bm{\gamma}_{\tau},
\end{equation}
where $Q_{\tau}(y|~\bm u,\bm x(t))$ is the $\tau$-th conditional quantile of $y$ given scalar covariates $\bm u$ and multiple functions $\bm x(t)$.
%$Q_{\tau}(y|~\bm u,\bm x(t))$ is the $\tau$-th conditional quantile of $y$ given covariates $\bm u$ and $\bm x(t)$,  $\alpha_\tau$ is the intercept, $\bm \gamma_{\tau}=(\gamma_{1\tau},\dots,\gamma_{p\tau})^T$ is a $p\times1$ vector of scalar coefficients, $\bm u(t)=(u_{1},\dots,u_{p})^T$ is a $p\times 1 $ vector of scalar covariates, $\bm \beta_{\tau}(t)=(\beta_{1\tau}(t),\dots,\beta_{m\tau}(t))^T$ is a $m\times 1$ vector of functional coefficients, and $\bm x(t)=\left(x_{1}(t),\dots,x_{m}(t)\right)^T$ is a $m\times 1$ vector of functional covariates. This framework is a natural generalization of the simple linear quantile regression model (\ref{quantile:model0}) and is capable of accommodating both multiple scalars and functionals in the same time.
To our best knowledge, only a few works have studied this model; for example, %to approximate the functional coefficients $\bm \beta_\tau(t)$,
\cite{yu2016partial} used partial quantile basis while \cite{yao2017regularized} used penalized principal component basis.
Inspired by the success of wavelet basis with regularization in functional linear model (\citealp{zhao2012wavelet}; \citealp{wang2014regularized}; \citealp{zhao2015wavelet}),
we use it to approximate the functional coefficients $\bm \beta_\tau(t)$ in model (\ref{quantile:model}).
Wavelet basis can provide a good representation of functional coefficients by using only a small number of basis and
are particularly useful for capturing localized functional features.
Moreover, the wavelet transform is computationally efficient and hence suitable for dealing with multiple functional predictors.

The penalization we impose is sparse group lasso (\citealp{Zhao2014quantile}, \citealp{simon2013sparse}), which is motivated by the attention deficit hyperactivity disorder (ADHD) study from the ADHD-200 Sample Initiative Project.
Our goal is to predict ADHD index at various quantile levels by using both demographic information and functional magnetic resonance imaging (fMRI) data,
where the fMRI data consists of $116$ functional features, each of which represents a single region of interests (ROI) of human brain.
The sparse group lasso technique, by imposing a convex combination of lasso and group lasso penalties,
can select important ROIs while capture shared information among them.
More specifically, the group lasso penalty makes a sparse selection out of $116$ functional features of ROIs,
while the lasso penalty induces a sparse representation of each feature.
Common wavelet basis is used to represent different features so that the shared information among them can be captured.

There are five major contributions of this paper.
First, our conditional quantile framework provides a more suitable modelling of reality especially %the relationship between response and predictors
when the response is heavy tailed (\citealp{yao2017regularized}).
It is also a compelling choice of dealing with heteroscedasticity issues and can provide a more complete picture of the response (\citealp{koenker2005quantile}).
Second, the wavelet basis we adopt provides a good approximation of functional coefficients while effectively detects the local features. % as well.
The wavelet transform we use is computationally efficient and hence can be easily extended to deal with multiple functional predictors.
Third, the proposed sparse group lasso method selects important functional predictors and retains shared information among them as well.
It is extremely useful in ADHD-200 fMRI study so that both individual and common information can be captured among the different ROIs.
Fourth, the estimation problem is in fact a penalized quantile regression problem,
which can be reformulated into a second-order cone program and then easily solved by an interior point method implemented by a powerful R package: Rmosek.
We also propose a novel algorithm to solve it by using alternating direction method of multipliers (ADMM).
Fifth, we successfully derive the asymptotic properties including the convergence rate and prediction error bound
which theoretically warrants good performance of our estimates.

The rest of paper is organized as follows.
In Section 2, we review some necessary background on wavelets and provide the penalized quantile objective function with sparse group lasso penalty.
The asymptotic properties such as the convergence rate and predictor error bound are established in Section 3.
In Section 4, the quantile penalization problem is reformulated into a second-order cone program (SOCP) and solved by an interior point method by using a powerful R package: Rmosek.
We also propose a novel algorithm using alternating direction method of multipliers (ADMM).
Finite sample simulations and a real data from ADHD-200 fMRI data are investigated in Section 5 to illustrate
the superiority of our proposed method.

\section{Wavelet-based Sparse Group Lasso}
In this section, we first review some necessary background on wavelets.
We then provide the penalized quantile objective function with sparse group lasso penalty
where the functional coefficients are approximated by wavelet basis.
This leads to the sparsities of both the selection and representation of functional features.
More specifically,
the group lasso selects a sparse set from available functional features,
while the lasso induces a sparse representation of the selected functional features.
% and then project our data into the space generated by the wavelet basis. After that, we take advantage of the sparse representation of the functions in the wavelet domain to derive our penalized objective function.

\subsection{Some Background on Wavelets}
Wavelets are basis function that can provide a good approximation of functional coefficients while effectively capture the local features (\citealp{zhao2012wavelet}).
Moreover, the wavelet transform is computationally efficient and hence can be easily extended to deal with multiple functional predictors (\citealp{daubechies1990wavelet}).
For a given $\tau \in (0,1)$, let $\beta_{l\tau}(t)$ be one component of $\bm \beta_\tau(t)$ in (\ref{quantile:model}), where
$\bm \beta_\tau(t) = ( \beta_{1\tau}(t), \ldots, \beta_{m\tau}(t) )^T$.
Suppose that $\beta_{l\tau}(t)$ is in $L^2[0,1]$. We can approximate it using wavelet basis.
For any wavelet basis in $L^2[0,1]$, they can be derived by dilating and translating two orthonormal basic functions: a scaling function and a wavelet function, namely $\phi(t)$ and $\psi(t)$ respectively:
%In particular,  for two orthonormal basic functions in ${L}^2[0 ,1]$,  let $ \varphi(t)$ be a scaling function and $\psi(t)$ be a wavelet function such that $\int_{0}^{1}\varphi(t)=1$ and  $\int_{0}^{1}\psi(t)=0$.
%The wavelet basis are obtained by dilating and translating $\varphi(t)$ and $\psi(t)$ respectively:
\begin{equation*}
	\varphi_{jk}(t)=\sqrt{2^{j}}\varphi(2^{j}t-k),\ \ \ \
	\psi_{jk}(t)=\sqrt{2^{j}}\psi(2^{j}t-k),
\end{equation*}
where $j$ and $k$ are integers, $\int_{0}^{1}\varphi(t)=1$ and  $\int_{0}^{1}\psi(t)=0$.
%The construction of a wavelet basis for ${L}_2[0 ,1]$ starts with two orthonormal basic functions: a scaling function, $ \varphi(t)$, and a wavelet function, $\psi(t)$, satisfying $\int_{0}^{1}\varphi(t)=1$ and  $\int_{0}^{1}\psi(t)=0$.
In particular, given a primary resolution level $j_0$, the wavelet basis are %defined as an orthonormal basis of ${L}^2[0 ,1]$:
\begin{equation}\label{wb}
	 \{ \varphi_{j_0,k}\}_{ 0\leq k\leq 2^{j_0}-1} \quad \rm{and} \quad \{\psi_{j,k}\}_{j_0 \leq j,\ 0\leq k\leq 2^{j}-1} .
\end{equation}
Therefore, $\beta_{l\tau}(t)$ can be approximated by
\begin{equation}\label{expansion1}
	\beta_{l \tau}(t)=\sum_{k=0}^{2^{j_0}-1}a_{j_0k}^l\varphi_{j_0k}(t)+\sum_{j=j_0}^{\infty
	}\sum_{k=0}^{2^j-1}d_{jk}^l\psi_{jk}(t),  \ \ \ \ \ \rm{for}~~ l=1,\dots,m,
\end{equation}
where $a_{j_0k}^l=\int_{0}^{1}\beta_{l\tau}(t)\varphi_{j_0,k}(t)dt$ is the approximation coefficients at the coarsest resolution $j_0$,
and $d_{jk}^l=\int_{0}^{1}\beta_{l\tau}(t)\psi_{jk}(t)dt$ is the detail coefficients characterizing the fine structures.

In practice, the functional covariates $\bm x(t) = ( x_1(t),\ldots,x_m(t))^T$ are discretely observed,
for instance without loss of generality, at $N=2^J$ equally spaced points of $[0,1]$ with
$0=t_1 < t_2 < \cdots < t_N = 1$.
Let $\bm X = (\bm x_1, \ldots,\bm x_m)$ and $\bm\beta_\tau = (\bm\beta_{1\tau}, \ldots, \bm \beta_{m\tau})$,
where $\bm x_l = ( x_1(t_1),\ldots,x_m(t_N))^T$, $\bm \beta_{l \tau} = ( \beta_{l \tau}(t_1),\ldots, \beta_{l \tau}(t_N))^T$ and $l = 1,\ldots,m$.
We represent $\bm X$ and $\bm \beta_\tau$ by the wavelet coefficients through discrete wavelet transform (DWT).
In particular, let $\bm W$ be an $N\times N$ matrix associated with orthonormal wavelet basis derived from DWT.
Suppose $\bm C$ and $\bm B$ are the corresponding wavelet coefficients of $\bm X$ and $\bm\beta_\tau$.
Then we have $\bm X = \bm W^T \bm C$, $\bm \beta_\tau = \bm W^T \bm B_\tau$,
and the integration in model (\ref{quantile:model}):
\begin{eqnarray*}
\int_{0}^{1}\bm x^T(t)\bm \beta_{\tau}(t)dt \approx \rm{vec}(\bm X)^T \rm{vec}(\bm \beta_\tau) /N = \rm{vec}(\bm W^T \bm C)^T \rm{vec}(\bm W^T \bm B) /N
= \rm{vec}(\bm C)^T \rm{vec}(\bm B)/N.
\end{eqnarray*}
The last equality holds due to the orthonormality of $\mathbf{W}$.
From now on, we denote $\bm v =\rm{vec}(\bm C)^T/N $ and $\bm \theta_\tau = \rm{vec}(\bm B)$
where $\bm C = (\bm c_1,\ldots,\bm c_m)$ and $\bm B = (\bm b_{1 \tau}, \ldots, \bm b_{m \tau})$.

\subsection{Model Estimation}
Using wavelet basis by DWT, model (\ref{quantile:model}) becomes
\begin{equation}\label{quantile:discrete}
	Q_{\tau}(y|\bm u,\bm x(t)) \approx \alpha_{\tau}+ \bm v^T \bm \theta_\tau + \bm{u}^T\bm{\gamma}_{\tau}.
	%+\varepsilon_i^*
\end{equation}
Given $n$ identical copies of data triplets $(\bm X_i, \bm u_i, y_i)$, where $\bm X_i$ and $\bm u_i$ are the observed functional and scalar covariates respectively, and $y_i$ is the corresponding response,
%A discrete approximation of model (\ref{quantile:model}) is
%\begin{equation}\label{quantile:discrete}
%	Q_{\tau}(y_i|\bm u_i,\bm x_i(t)) \approx \alpha_{\tau}+ \bm C_i^T \bm \Theta_\tau + \bm{u}_i^T\bm{\gamma}_{\tau},
%	%+\varepsilon_i^*
%	\ \ \ \rm{for}~~ i=1,\dots, n.
%\end{equation}
%For the notation simplicity, we hereafter denote $\bm \theta_{l\tau} = \bm \beta_{l\tau}^*$ for $l=1,\ldots,m$.
the parameters in (\ref{quantile:discrete}) can be estimated by minimizing a regular quantile loss function.
However, to find the important functional covariates in predicting responses while preserve a desired sparse representation of the coefficients,
an appropriate penalty has to be imposed.
In this paper, we propose to use the sparse group lasso penalty
\begin{equation}
\label{sglp}
P_{\lambda_1, \lambda_2}(\bm \theta) = \lambda_1\sum_{l=1}^{m}||\bm b_l||_1 +\lambda_2\sum_{l=1}^{m}||\bm b_l ||_2,
\end{equation}
where $||\cdot||_1$ and $||\cdot||_2$ represent the $L_1$ and $L_2$ norms respectively, and $\lambda_1$ and $\lambda_2$ are two nonnegative tuning parameters.
The sparse group lasso penalty includes two components, namely a lasso and a group lasso penalties,
where the lasso penalty $||\cdot||_1$ induces sparsity in each functional coefficient and
the group lasso penalty $||\cdot ||_2$ selects functional coefficients.
%The group lasso penalty $||\bm b_l ||_2$ selects a sparse set out of $\{ \bm b_1,\ldots,\bm
% b_m\}$ while the lasso penalty $||\bm b_l||_1$ induces sparsity in each $\bm b_l$.
Common information among functional covariates can be retained by using the same wavelet basis to approximate the functional coefficients.
Moreover, the sparse group lasso warrants the selection of important functional coefficients while captures distinct traits carried by individual functional covariates.
Specifically, the parameters $\alpha_{\tau}$, $\bm \gamma_{\tau}$, and $\bm \theta_{\tau}$ can be estimated by %minimizing a penalized quantile loss function.
%That is
%The group lasso part achieves selection of the functional predictors based on the preset grouping structure of the parameters $\bm \theta_{j\tau}$s.
%The Lasso part imposes smoothness of the coefficients $\bm \theta_{j\tau}$s by controlling the number of basis functions. In fact, this sparsity is also implied by using the wavelet basis. Thus, the parameters $\alpha_{\tau}$, $\bm \gamma_{\tau}$, and $\bm \theta_{\tau}$ can be estimated by  minimizing the loss function of quantile regression with some shrinkage constraints. That is
\begin{equation}\label{qr:sgl}
	(\hat{\alpha}_{\tau},\hat{\bm \gamma}_{\tau},\hat{\bm \theta}_{\tau})= \underset{\alpha,\bm \gamma,\bm \theta}{\arg\min} \sum_{i=1}^{n}\rho_{\tau}\big( y_i-\alpha -\bm u_i^T\bm \gamma -\bm v_i^T \bm \theta \big) + P_{\lambda_1, \lambda_2}(\bm \theta),
\end{equation}
where $\rho_\tau(x) = x(\tau-\mathbf{1}(x <0))$ is the quantile check function (\citealp{koenker2005quantile}).

To combine information from different quantiles, \citet{zou2008composite} proposed composite quantile regression,
%In particular, under the condition that the effects of covariates are piecewise constant or continuous across different quantile levels,
%composite quantile regression
which simultaneously considers multiple regression quantiles at different levels.
With homoscedasticity assumption, where all conditional regression quantiles have the same slope,
the composite quantile estimate is more efficient than the one from a single level and
has in recent years begun to gain its popularity in many fields (\citealp{kai2010local}; \citealp{fan2010selective}; \citealp{bradic2011penalized}, \citealp{kai2011new};  \citealp{yu2016partial}).
%In particular, under the homoscedasticity assumption, all conditional regression quantiles are parallel and we have the same coefficients but different intercepts.
In this paper, we propose to use composite quantile regression with sparse group lasso penalty in our functional data analysis framework.
Let $0<\tau_1<\dots<\tau_k<1$ denote the selected quantile levels and then
the parameters ${\bm{\alpha}},{\bm\gamma}$ and ${\bm{\theta}}$ can be estimated by
\begin{equation}\label{cqr:sgl}
	(\hat{\bm\alpha},\hat{\bm{\gamma}},\hat{\bm\theta})=\underset{{\bm{\alpha}},{\bm\gamma},{\bm{\theta}}} {\arg\min}\sum_{k=1}^{K} \sum_{i=1}^{n} {\rho}_{\tau_k}(y_i-\alpha_{k}-\bm{u}_i^T\bm{\gamma}
-\bm v_i^T \bm \theta)+P_{\lambda_1, \lambda_2}(\bm\theta),
\end{equation}
where $\bm\alpha=(\alpha_{1}.\dots,\alpha_{{K}})$ is a vector of intercepts.
Typically, we can choose $K = 9$ and use equally spaced quantiles (\citealp{kai2010local}; \citealp{zou2008composite}).
%For notation simplicity, we can denote	$L_n({\bm{\alpha},{\bm{\theta}},{\bm{\gamma}}})=\sum_{k=1}^{K} \sum_{i=1}^{n} {\rho}_{\tau_k}(y_i-\alpha_k-\bm C_i^T \bm \theta-\bm{u}_i^T\bm{\gamma})$ to be the objective function of composite quantile regression. %, which is a mixture of the objective functions from different quantile regression models.
Note that quantile estimate (\ref{qr:sgl}) at a single level is just a special case of composite quantile estimate (\ref{cqr:sgl}) with $K=1$.
%Therefore, we hereafter only focus on the composite quantile regression cases.
In the following, we will focus on the composite quantile regression case of (\ref{cqr:sgl}).

\section{Asymptotics}
In this section, we investigate the asymptotic properties of our proposed estimates when both the sample size $n$ and the number of discrete points $N_n$ tend to infinity.
% meaning that the sample size $n$ increases and the curves, $\bm x_i(t)$s, are also more densely observed, respectively.
Let $\lambda_{1,n}$ and $\lambda_{2,n}$ denote the tuning parameters %$\lambda_1$ and $\lambda_2$ respectively
when the sample size is $n$.
To derive the asymptotic properties, we impose the following conditions:
\begin{enumerate}
	\item[] {\bf A1.} \ \ The model errors $\varepsilon_1,\ldots,\varepsilon_n$ are independently following a distribution $ F$,
    with density $ f$ to be bounded away from zero and infinity,
    and its derivative $f^\prime$ to be continuous and uniformly bounded.
	\item[] {\bf A2.} \  \ There exist two constants $c_1$ and $c_2$ such that % the regression matrix $\bm A$ satisfies the eigenvalue condition
	\begin{equation*}
		0<c_1<\varrho_{\min}(\frac{1}{n}\bm{A}_n^T\bm A_n)\leq \varrho_{\max}(\frac{1}{n}\bm{A}_n^T\bm A_n)<c_2<\infty,
	\end{equation*}
	where $\bm A_n = (\bm a_1, \ldots, \bm a_n )^T$ is the design matrix with $\bm a_i=(1,\bm v_i^T,\bm u_i^T)^T$,
and $\varrho_{\min}(\cdot)$ and $\varrho_{\max}(\cdot)$ are the smallest and largest eigenvalues of $\frac{1}{n}\bm{A}_n^T\bm A_n$ respectively.
	\item[] {\bf A3.} \  \ There exists a constant $M$ such that $||\bm a_i||_{2}<M$ for all $i$.
	\item[] {\bf A4.} \  \ The functional slope $\beta_l(t)$s are $d$ times differentiable in the Sobolev sense, and the wavelet basis has $w$ vanishing moments, where $w>d.$
	\item[] {\bf A5.} \ \ $\lambda_{1,n}=O(\sqrt{n})$ and \ $\lambda_{2,n}=O(\sqrt{n})$. % and $n=O_p(N^{4d}_n)$.
	\item[] {\bf A6.} \ \ $N_n / n\to 0$.
\end{enumerate}
These regularity conditions might not be the weakest ones but are commonly assumed among literatures of quantile regression and functional linear model. % with possibly growing number of parameters.
Condition (A1) is standard for quantile regression (\citealp{koenker2005quantile}; \citealp{Zhao2014quantile}), which regulates the behavior of the conditional density of the response in a neighborhood of the conditional quantile and is crucial to the asymptotic properties of quantile estimators (\citealp{koenker1978regression}). Condition (A2) is a classical condition in functional linear regression literature (\citealp{delaigle2012methodology}). It ensures the eigenvalues of the covariance matrix go to neither zero nor infinity too quickly. %not too quickly go to either zero or infinity.
Similar conditions as (A3) - (A6) can be found in \cite{zhao2012wavelet} and \cite{zhao2015wavelet}, among others.
Condition (A4) guarantees that the space spanned by the wavelet basis can well approximate the functional slopes with small approximation errors. % go to 0.
Condition (A6) implies that to allow for estimation of $\beta$ with appropriate asymptotic properties, $n$ should grow faster than $N_n$.
%In other words, as both $n$ and $N_n$ grow large, i.e. as both the sample size and the resolution at which the predictors are observed increase,
%the improvement of resolution requires sufficient ``more information'' to have a good estimation of $\beta$.
Note the wavelet basis has $w$ vanishing moments if and only if its scaling function $\varphi$ can generate polynomials of degree at most $w$.
\begin{thm}\label{cqr:bound}
	Let $\hat{ \beta}_{l,n}$ be the estimator resulting from (\ref{cqr:sgl}) and $\beta_l$ is the true coefficient function. If Conditions (A1)-(A6) hold, then
	\begin{equation*}
		||\hat{ \beta}_{l,n}-\beta_l||_2^2=O_p\left({\frac {N_n} n}\right)+o_p\left(\frac 1 {N_n^{2d}}\right).
	\end{equation*}
\end{thm}
A detailed proof of this theorem is provided in the Appendix.
The accuracy of $\hat{\beta}$ relies on both $n$ and $N_n$.
The approximation error rate of $\hat{\beta}$ towards $\beta$ are controlled by two terms.
The first term is of the same order of $N_n / n$ which is a typical result of estimating,
while the second term is of the lower order of $1/N_n^{2d}$ which is mainly due to approximation by wavelets.
In particular, the approximation error rate is dominated by the second term if $N_n^{2d+1}$ is of the lower order of $n$.
Otherwise, it is dominated by the first term.
Under some further conditions, we can have the following theorem for the prediction error bound:
\begin{thm}% For a given fixed $k$,
	Suppose $ x_l(t)$ is square integrable on $[0,1]$ and $\bm F^{-1}(\tau)=0$. If Conditions (A1)-(A6) hold and $F^{-1}(\tau_k)=0$, then
	\begin{equation*}
		||\hat{y}-y||_2^2=O_p({\frac {N_n} n} )+o_p(\frac 1 {N_n^{2d}}),
	\end{equation*}
	where $y$ is the true response and $\hat{y}$ is estimated $\tau_k$'s conditional quantile.  %from (\ref{cqr:sgl}).
\end{thm}
The proof follows that from Theorem \ref{cqr:bound} and the Cauchy-Schwarz inequality, the details of which are omitted in this paper.
Similarly as in Theorem \ref{cqr:bound},
%Note the quantile effect of composite quantile prediction is only reflected by the quantile intercept $\alpha_k$.
%Although the true response $y$ may depend on the functional covariates $\mathbf{x}(t)$ through model (\ref{quantile:model}),
$L_2$ prediction error rate depends on the same two terms from estimating and approximation by wavelets respectively,
while the estimation errors caused by $\hat{\alpha}_k$ and $\hat{\gamma}$ is absorbed by the first term.
%the sample size $n$ and the number of discrete points $N_n$, and is controlled by two terms which are the same as for Theorem \ref{cqr:bound}.
%The first term results from estimating, the errors of which are absorbed by $\hat{\alpha}_k$ and $\hat{\gamma}$,
%while the second term is mainly due to approximation by wavelets.
%Note the prediction error rate is controlled by the second term if $N_n$ is of the lower order of $n^{\frac{1}{2d+1}}$.
%Otherwise, it is controlled by the first term.
%--------------------------------------  Proof, don?t delect, move the appendix ---------------------------------------------%
%\begin{proof}
%	The proof follows from \ref{cqr:bound} and the Cauchy-Schwarz inequality. We omit the details.
%\end{proof}	
%--------------------------------------  Proof, don?t delect, move the appendix ---------------------------------------------%

\section{Implementations}
Due to the non-smoothness of loss function, quantile estimator does not enjoy the nice asymptotic properties, as well as computational easiness, as what ordinary least square estimator does.
After illustrating asymptotic theory of the proposed quantile estimator, it becomes of great importance to have an efficient algorithm to obtain it.
In this section, we reformulate the optimization problem (\ref{cqr:sgl}) into a second-order cone program (SOCP)
and implement it by interior point method using a powerful R package: {\bf Rmosek} (\citealp{aps2015rmosek}).
Alternatively we propose a novel algorithm to solve problem (\ref{cqr:sgl}) by using alternating direction method of multipliers (ADMM)
which was a technique recently employed by many researchers in solving penalized quantile regression problems.
In the end, we discuss some practical rules to choose tuning parameters. %$\lambda_1$ and $\lambda_2$.

\subsection{A Second-Order Cone Program}
%\subsection{An Interior Point Algorithm}
%The optimization problem (\ref{cqr:sgl}) can be reformulated into a standard form of second-order cone programing (SOCP).
Let the superscripts ${}^{+}$ and  ${}^{-}$ denote the positive and negative parts of a vector.
%Let $\bm\theta^+  $ and $\bm\theta^-$ be the positive and negative parts of $\bm\theta$,
%with $\bm b_l^+  $ and $\bm b_l^-$ be the positive and negative parts of $\bm b_l$.
For unknown parameter $\bm\theta$ in (\ref{cqr:sgl}), we write: $\bm\theta=\bm\theta^+-\bm\theta^-$ and $||\bm\theta||_1=||\bm\theta^+||_1+||\bm\theta^-||_1$.
Similarly, we have $\bm b=\bm b^+-\bm b^-$ and $||\bm b||_1=||\bm b^+||_1+||\bm b^-||_1$.
Then problem (\ref{cqr:sgl}) can be reformulated as the following standard second-order cone program:
\begin{eqnarray}{\label{SOCP3}}
	\min&&\sum_{k=1}^K\sum_{i=1}^{n}(\tau_k r_{ki}^+ +(1-\tau_k) r_{ki}^-)+\lambda_1\sum_{l=1}^{m}(||\bm b^+_l||_1+||\bm b^-_l||_1)+\lambda_2\sum_{l=1}^{m}z_l\nonumber \\
	\textrm{subject to}&&    -r_{ki}^-\leq y_i-\alpha_k-\bm{u}_i^T\bm{\gamma}-\bm{v}_i^T(\bm\theta^+ - \bm\theta^-)\leq r_{ki}^+\nonumber \\
	&&      \sqrt{||\bm b^+_{l}||_2^2+||\bm b_{l}^-||_2^2}\leq z_l\\
	&& \bm\theta^+\succeq 0,~ \bm\theta^-\succeq  0,~ z_l\geq 0,~ r_{ki}^+\geq 0,\ r_{ki}^-\geq 0.\nonumber
\end{eqnarray}
where $r_{ki}^+$, $r_{ki}^-$ and $z_l$ are three nonnegative slack variables,
%such that
%$\bm\theta=\bm\theta^+-\bm\theta^-$ and $||\bm\theta||_1=||\bm\theta^+||_1+||\bm\theta^-||_1$.
%Each term in the objective function of composite quantile regression can be written as a linear constraint; for instance,
%\begin{equation*}
%	-r_{ki}^-\leq y_i-\alpha_k-\bm{u}_i^T\bm{\gamma}-\bm{v}_i^T(\bm\theta^+-\bm\theta^-)\leq r_{ki}^+, \ \ \ \ 	 \\
%\end{equation*}
 %Finally, we deal with the $L_2$ type penalty in $P_{\lambda_1,\lambda_2}(\bm\theta)$.
and the contraint $\sqrt{||\bm b^+_{l}||_2^2+||\bm b_{l}^-||_2^2}\leq z_l$ implies
a second order cone of dimension $2N+1$ (\citealp{lobo1998applications}) denoted as
\begin{equation*}
	\mathbb{Q}^{2N+1}_l=\bigg\{(z_l, \bm b^+_{l},\bm b^-_{l})\in \mathbb{R}^{2N+1} \big | \ z_l\geq \sqrt{||\bm b^+_{l}||_2^2+||\bm b_{l}^-||_2^2} \bigg\}.\
\end{equation*}
%Once we have the solutions of $\bm\theta^+$ and $\bm\theta^-$, we can cosntruct $\bm\theta$ by $\bm\theta=\bm\theta^+-\bm\theta^-$.	
The reformulation is guaranteed by the fact that for each component of optimal $\bm b_l$, either $ b^+_{l,j}=0$ or  $ b^-_{l,j}=0$
would be held. Otherwise, for optimal $\bm b_l$, if there exist $l$ and $j_0$ such that $b^+_{l,{j_0}}>0$ and  $b^-_{l,{j_0}}>0$, we can replace $b^+_{l,{j_0}}$ and  $b^-_{l,{j_0}}$ by $b^{\rm{(new)}+}_{l,{j_0}}$ and  $b^{\rm{(new)}-}_{l,{j_0}}$ respectively with
\begin{equation*}
	b^{\rm{(new)}+}_{l,{j_0}}=\left\{
	\begin{array}{ll}
		0\ &if\ b^+_{l,{j_0}}<b^-_{l,{j_0}},\\
		b^+_{l,{j_0}}-b^-_{l,{j_0}} &  \rm{otherwise},
	\end{array}
	\right.\ \ \
	b^{\rm{(new)}-}_{l,{j_0}}=\left\{
	\begin{array}{ll}
		0\ &if\ b^+_{l,{j_0}}>b^-_{l,{j_0}},\\
		b^-_{l,{j_0}}-b^+_{l,{j_0}}, &  \rm{otherwise}.
	\end{array}
	\right.
\end{equation*}
As a result, the objective function in (\ref{SOCP3}) decreases, which contradicts with the fact that $\bm b_l$ being optimal.
%Therefore, we have $||\bm b_l ||_1 = ||\bm b^+_l-\bm b^-_l||_1 = ||\bm b^+_l||_1+||\bm b^-_l||_1$ and  $||\bm b_l ||_2 = ||\bm b^+_l-\bm b^-_l||_2 = \sqrt{||\bm b_l^+||_2^2+||\bm b_l^-||_2^2}$.

%After reformulating the original penalized composite quantile regression problem (\ref{cqr:sgl}) into a SOCP (\ref{SOCP3}),
Various optimization strategies can be applied to solve SOCP (\ref{SOCP3}) such as interior point method (\citealp{koenker1996interior}) and the simplex method (\citealp{koenker2005quantile}).
In this paper, we choose to use interior point method.
The R package we use is \textbf{Rmosek} (\citealp{aps2015rmosek}).
The technique proposed to reformulate our problem into a SOCP can be easily adapted to other penalized quantile regression problems;
for example, quantile ridge regression (\citealp{wu2009variable}).
%, we can obtain the estimations by solving the equivalent form~(\ref{SOCP3}).

	\subsection{ADMM Algorithm}
	
	%\subsection{Overview of ADMM}
Although problem (\ref{cqr:sgl}) is convex,
%the sum of a composite quantile loss function and a sparse group lasso penalty, both of which are convex.
solving it can be very slow partially due to large scale data in the application
and the non-smooth terms in the objective that prevent fast gradient method being applied.
However, with non-smooth terms in the objective and very large scale data, these methods can be very slow.
In this section, we explore the additive structure of the objective function,
namely, decompose it into two sub convex problems, and then propose a novel and efficient algorithm by using alternating direction method of multipliers (ADMM) (\citealp{gabay1976dual}).
This powerful tool was originated in 1950s and developed during 1970s (\citealp{hestenes1969multiplier}; \citealp{gabay1976dual}).
It has been popularized in recent years among quantile regression literature (\citealp{boyd2011distributed}; \citealp{gao2015cqrreg}; \citealp{kong2015estimation}).

    Denote $L_n({\bm{\alpha},{\bm{\theta}},{\bm{\gamma}}}) = \sum_{k=1}^{K} \sum_{i=1}^{n} {\rho}_{\tau_k}(y_i-\alpha_{k}-\bm{u}_i^T\bm{\gamma} - \bm v_i^T \bm \theta)$.
	The minimization problem (\ref{cqr:sgl}) can be rewritten as
	\begin{eqnarray*}
	\text{min} & & L_n({\bm{\alpha},{\bm{\theta}},{\bm{\gamma}}})+P_{\lambda_1, \lambda_2}(\bm\theta^*)\\
	\text{subject to} & & {\bm{\theta}}={\bm{\theta^*}},
	\end{eqnarray*}
	where $L_n(\cdot)$  and $P_{\lambda_1, \lambda_2}(\cdot)$ are two convex functions.
	Applying augmented lagrangian \citep{hestenes1969multiplier}, we have
	\begin{equation}
	\begin{aligned}
	L_{n,\eta}({\bm{\alpha},{\bm{\theta}},{\bm{\gamma}},\bm\theta^*,\bm{\mu}})=& L_n({\bm{\alpha},{\bm{\theta}},{\bm{\gamma}}})+P_{\lambda_1, \lambda_2}(\bm\theta^*) + \bm \mu^T(\bm\theta - \bm\theta^*)
	+\frac{\eta}{2}\| \bm\theta - \bm\theta^*\|_2^2.
	\end{aligned}
	\label{eqnLargn}
	\end{equation}
%		The  ADMM algorithm to obtain the minimizer of (\ref{eqnLargn}) follows a three-step iterative scheme:
%	\begin{eqnarray*}
%	\label{eq:t}
%	 ({\bm{\alpha}^{(l+1)},{\bm{\theta}}^{(l+1)},{\bm{\gamma}}^{(l+1)}}) &=& \underset{\bm{\alpha},{\bm{\theta}},{\bm{\gamma}}}{\text{argmin}}  \ L_{n,\eta}({\bm{\alpha},{\bm{\theta}},{\bm{\gamma}},\bm\theta^{*(l)},\bm \mu^{(l)}})
%	\\
%	 \bm\theta^{*(l+1)} &=&\underset{\theta^*}{\text{ argmin}} \ L_{n,\eta}({\bm{\alpha}^{(l+1)},{\bm{\theta}}^{(l+1)},{\bm{\gamma}^{(l+1)}},\bm\theta^{*},\bm \mu^{(l)}})
%	\\
%	 \bm \mu^{(l+1)}&=&\bm \mu^{(l)}+ \eta( \bm \theta^{(l+1)} - \bm \theta^{*(l +1)}).
%	\end{eqnarray*}
 Let $ \bm w= \bm\mu / \eta.$
	The ADMM algorithm to obtain the minimizer of (\ref{eqnLargn}) follows a three-step iterative scheme:
	\begin{eqnarray}
	\label{eq:12}
	 ({\bm{\alpha}^{(l+1)},{\bm{\theta}}^{(l+1)},{\bm{\gamma}}^{(l+1)}}) \nonumber
	&=& \underset{\bm{\alpha},{\bm{\theta}},{\bm{\gamma}}}{\text{argmin}} \ L_n({\bm{\alpha},{\bm{\theta}},{\bm{\gamma}}})+\frac{\eta}{2}\|{\bm{\theta}}-{\bm{\theta}}^{*(l)}+\bm w^{(l)} \|_2^2\\
	\bm\theta^{*(l+1)} &=&\underset{\theta^*}{\text{ argmin}} \ P_{\lambda_1, \lambda_2}(\bm\theta^*)
	+\frac{\eta}{2}\|{\bm{\theta}}^{(l+1)}-{\bm{\theta}}^{*}+\bm w^{(l)} \|_2^2 \nonumber \\
	\bm w^{(l+1)} &=& \bm w^{(l)}+\eta({\bm{\theta}}^{(l+1)}-{\bm{\theta}}^{*(l+1)}).
	\end{eqnarray}
		For the first step of (\ref{eq:12}), it can be reformulated as a SOCP:
	\begin{eqnarray*}
	 %\underset{\bm{\alpha},{\bm{\theta}},{\bm{\gamma}}}
{\text{min}} &&\sum_{k=1}^{K}\sum_{i=1}^{n}\rho_{\tau_k}(r_{ik})+\frac{\eta}{2} \|{\bm{\theta}}-{\bm{\theta}}^{*(l)}+\bm w^{(l)} \|_2^2  \\
	\text{subject to} && {y_i-\alpha_k-\bm u_i^T\bm \gamma -\bm v_i^T \bm \theta =r_{ik}},\ \ \rm{for}
	\ i=1,
	 \dots, n;\ \ k=1,\ \dots,\ K,
	\end{eqnarray*}
	which can be easily solved by following an ADMM scheme:
		\begin{eqnarray}
		\label{eq:u}
		r_{ik}^{(j+1)}&=&\underset{r_{ik}}{\text{ argmin }}\rho_{\tau_k}(r_{ik})+\frac{\eta_1}{2}(y_i-\alpha_k^{(j)}-\bm u_i^T\bm \gamma^{(j)} -\bm v_i^T \bm \theta^{(j)}-r_{ik}+ z  ^{(j)}_{ik})^2 \nonumber \\
		({\bm{\alpha}^{(j+1)},{\bm{\theta}}^{(j+1)},{\bm{\gamma}}^{(j+1)}})&=&\underset{\bm{\alpha},{\bm{\theta}},{\bm{\gamma}}}{\text{argmin}} \frac{\eta}{2} \|{\bm{\theta}}-{\bm{\theta}}^{*(l)}+\bm w^{(l)} \|_2^2+\frac{\eta_1}{2}\sum_{k=1}^{K}\sum_{i=1}^{n}(y_i-\alpha_k-\bm u_i^T\bm \gamma -\bm v_i^T \bm \theta-r_{ik}^{(j+1)}+z^{(j)}_{ik})^2 \nonumber \\
		z^{(j+1)}_{ik}&=& z^{(j)}_{ik}+\eta_1(y_i-\alpha_k^{(j+1)}-\bm u_i^T\bm \gamma^{(j+1)} -\bm v_i^T \bm \theta^{(j+1)}-z_{ik}^{(j+1)}).
		\end{eqnarray}
	   The first step of  (\ref{eq:u}) can be explicitly solved by the soft thresholding operator.
	   The second step can be easily approximated by a standard ridge regression therefore has a closed form.
	
	The second step of (\ref{eq:12}) can be simplified by the soft thresholding operator. That is,
	\begin{eqnarray*}
	%\label{eq:12}
	\bm \upsilon^{*} &=&\ sgn({\bm{\theta}}^{(l+1)}+\bm w^{(l)}) \cdot \max(|\bm{\theta}^{(l+1)}+\bm w^{(l)}|-\frac{\lambda_1}{\eta},\ 0)\\
\bm\theta^{*(l+1)} &=&\ \frac{\bm \upsilon^{*}}{||\bm \upsilon^{*}||_2}\max(||\bm \upsilon^{*}||_2-\frac{\lambda_2}{\eta}, 0),
	\label{eqnThresh}
	\end{eqnarray*}
	where $sgn(\cdot)$ is the sign function.

%		Therefore, the ADMM algorithm of (\ref{eq:12}) can be much faster and more efficient than the other techniques.

	A typical stopping criterion with primal and dual residuals denoted respectively by $r_{\rm{primal}}$ and $r_{\rm{dual}}$ (\citealp{boyd2011distributed}) can be chosen as :
	$$||\bm{\theta}^{(l)}-\bm{\theta}^{*(l)}||_2\leq
	r_{\rm{primal}} \,\,\,\,\,\text{and  }\,\,\,\,\|\eta(\bm{\theta}^{*(l)}-\bm{\theta}^{*(l-1)})   \|_2\leq r_{\rm{dual}},$$
	with
	\begin{eqnarray*}
	 r_{\rm{primal}} &=& \sqrt{mN} \epsilon_{\rm{abs}}+\epsilon_{\rm{rel}} \cdot \max\, \{\|\bm{\theta}^{(l)}\|_2,\,\|\bm{\theta}^{*(l)}\|_2\},
	\\
	 r_{\rm{dual}}&=& \sqrt{mN+q+K}\epsilon_{\rm{abs}}+\epsilon_{\rm{rel}} \cdot \|\bm w^{(l)}\|_2,
	\end{eqnarray*}
	where $q$ is the dimension of $\bm\gamma$, and parameters $\epsilon_{\rm{abs}}$ and $\epsilon_{\rm{rel}}$ are two predefined absolute and relative tolerances which can be set as $10^{-4}$ and  $10^{-2}$ respectively.

%    The advantage of using ADMM is that
    Instead of tackling the original problem directly, ADMM decompose it into several sub convex problems then deal with them separately by iteration.
    In each iteration, the sub problem can be easily and efficiently solved by the soft thresholding operator or approximated to have a closed form.
    Therefore, the ADMM algorithm derived is much faster and more efficient than other general techniques.

\subsection{Selection of Tuning Parameters}
The proposed method involves selection of two nonnegative tuning parameters, namely $\lambda_1$ and $\lambda_2$, which control the severity of penalization towards model complexity.
%the degree of sparsity  individual and group sparsities of model.
Specifically, $\lambda_1$ controls sparsity in each functional coefficient while $\lambda_2$ controls the number of selected functional coefficients.
%For instance, if we expect strong overall sparsity and would like to encourage selecting coefficients, we can set relative large $\lambda_1$ and small $\lambda_2$.
%On the other hand, if we expect only a small number of selected coefficients, with a mild sparsity within each of them, we can use small $\lambda_1$ and large $\lambda_2$.
Although many options exist for selecting tuning parameters, such as AIC, BIC and cross validation,
there is no agreed-upon selection criterion in general.
%In practice, people usually pre-specify a finite set of values for the regularization parameters, then use either a validation dataset or a certain model selection criterion to pick the regularization parameters.
%Although some existing criteria such as AIC and cross validation could be potentially useful, \citet{wang2007tuning} and \citet{zhang2010regularization} showed that both of them may fail to consistently identify the true model. %hence become problematic.
After showing that AIC and cross validation may fail to consistently identify the true model,
\citet{zhang2010regularization} proposed to use the generalized information criterion (GIC), encompassing the commonly used AIC and BIC,
and illustrated the corresponding asymptotic consistency. %, to select tuning parameters.
More recently, \citet{zheng2015globally} used the GIC to make consistent model selection for quantile regression in ultra-high dimensional settings.
In this paper, we propose to use the GIC:
\begin{equation}
\label{gic}
	(\hat{\lambda}_1,\hat{\lambda}_2)=\underset{\lambda_1, \lambda_2}{\arg\min}\frac{1}{K}
	\sum_{k=1}^{K} \ln \left (\frac 1 n\sum_{i=1}^{n}\rho_{\tau_k}\left(y_i-\hat{y}_{ki}\right)\right) + \phi_n||\hat{\bm\theta}_{{\lambda}_1,{\lambda}_2}||_0,
\end{equation}
where $\hat{\bm \theta}_{{\lambda}_1,{\lambda}_2}$ is a solution of problem (\ref{cqr:sgl}), $||\cdot||_0$ denotes $L_0$ norm (total number of non-zero elements in a vector), $\phi_n$ is a sequence converging to zero with $n$ goes to infinity, and  $\hat{y}_{ki}$ is calculated from  $(\ref{quantile:discrete})$ with $\tau=\tau_k$.

In addition, we can also use the validation set  (\citealp{li2007quantile}, \citealp{wu2009variable}) to select gold standard $\lambda_1$ and $\lambda_2$ that minimize the prediction error. Simulations in Section 5 demonstrate a satisfactory behavior of the proposed criterion compared with the validation set method. % for selecting tuning parameters.

\section{Numerical Studies}
In this section, we compare performances of the proposed sparse group lasso method with group lasso and lasso methods using simulations and a real data from ADHD-200 fMRI sample (\citealp{mennes2013making}).
%The merits we use include coefficient estimation error and mean absolute prediction error.
We also compare the tuning parameters selected by the GIC approach we proposed and the validation set approach.
In our numerical studies,
we employ least-asymmetric wavelets of Daubechies with 6 vanishing moments and fix the tuning parameter ratio ${\lambda_1}/{\lambda_2}=0.5$ (\citealp{simon2013sparse}).
To simplify notations, we use qSGL, qL and qGL to represent the quantile sparse group lasso, lasso and group lasso methods respectively.

\subsection{Simulations}
%In this section, we compare the finite sample performance of a number of different methods, quantile Lasso and quantile group lasso, with regard to the functional estimation error and the  prediction error. We also examine the effectiveness of the proposed modeling strategy by investigating whether our method selects functional predictors appropriately. We use least-asymmetric wavelet of Daubechies with 6 vanishing moments for both the simulation study and the real data analysis, and fix the ratio ${\lambda_1}/{\lambda_2}=0.5$.

Our data are randomly generated using $12$ functional covariates and $2$ scalar covariates in a setting similar to \citet{collazos2016consistent}.
%The simulations are based on $200$ and $400$ generated observations of 12 functional covariates and 2 scalar covariates with a scalar response, extending the simulation setup of \citet{collazos2016consistent} by including more functional predictors.
In particular, the model is of the form:
\begin{eqnarray*}
	y_i&=&\alpha+\bm{u}_i^T\bm{\gamma}+	\int_{0}^{1}\bm x_{i}(t)^T \bm \beta(t)dt+\sigma\varepsilon_i\ \ \ \textrm{for}\ i=1 \dots, n,
\end{eqnarray*}
where $\bm u_i=(u_{i1},u_{i2})^T$ with $u_{i1}\sim N(0,1)$ and $u_{i2}\sim \textrm{Bernoulli}\left(0.5\right)$, and the coefficients $\bm{\gamma}=(0.32/256,\ 0.32/256)^T$.  The functional covariates $\bm x_i(t) = ( x_{i1}(t),\ldots,x_{i 12}(t))^T$ are observed on an equally spaced grid of $N=256$ points on $[0,1]$ with
\begin{eqnarray*}
	x_{i1}(t)&=&\sqrt{.84} \omega_{i1}(t)+ .4 \omega_{i6}(t),\ \ \  	x_{i2}(t)=\sqrt{.98} \omega_{i2}(t)+ .1 \omega_{i1}(t)+.1  \omega_{i5}(t),\\
	x_{i3}(t)&=&\sqrt{.84} \omega_{i3}(t)+ .4  \omega_{i4}(t),\ \ \   x_{i5}(t)=\sqrt{.99} \omega_{i5}(t)+ .1 \omega_{i2}(t),\\
	x_{il}(t)&=&\omega_{il}(t)\ \ \ \ \ \textrm{for}\ \ l=4, 6, 7,\dots, 12;
\end{eqnarray*}
where
\begin{equation*}
	\omega_{il}(t)=z_{il}(t)+\epsilon_{il}, ~ \epsilon_{il}\sim N\left(0, \left(.05r_{x_{il}}\right)^2\right), \ \ \ \textrm{for}\ l=1 \dots, 12,
\end{equation*}
with $r_{x_{il}}=\underset{i}{\max} \left( z_{il}(t) \right)-\underset{i}{\min}\left( z_{il}(t)\right)$ and
\begin{eqnarray*}
	\ z_{i1}(t)&=&\cos (2\pi(t-a_1))+a_2, \mathbb{T}_1=[0,~1],~ a_1\sim N(-4,~ 3^2),~ a_2\sim N(7,~1.5^2),\\
	\ z_{i2}(t)&=&b_1t^3+b_2t^2+b_3t, \mathbb{T}_2=[-1,~1],~ b_1\sim N(-3,~ 1.2^2),~ b_2\sim N(2, .5^2),~ b_3\sim N(-2~,1),\\
	\ z_{i3}(t)&=&\sin(2(t-c_1))+c_2t, \mathbb{T}_3=[0,~\pi/3],~ c_1\sim N(-2,~ 1),~ c_2\sim N(3, 1.5^2),\\
	\ z_{i4}(t)&=&d_1\cos(2t)+d_2t, \mathbb{T}_4=[-2,~1],~ d_1\sim U(2,~ 7),~ d_2\sim N(2, .4^2),\\
	\ z_{i5}(t)&=&e_1\sin (\pi t)+e_2, \mathbb{T}_5=[0,~\pi/3],~ e_1\sim U(3,~ 7),~ e_2\sim N(0,1),\\
	\ z_{i6}(t)&=&f_1e^{-t/3}+f_2t+f_3, \mathbb{T}_6=[-1,~1],~ f_1\sim N(4,~ 2^2),~ f_2\sim N(-3, .5^2),~ f_3\sim N(1,1),\\
	\ z_{il}(t)&=& 5\sqrt{2} \sum\nolimits_{j=1}^{49}\cos(j\pi t)g_j+5 h, \mathbb{T}_l=[0,~1],~g_j\sim N\left(0,~ \left({j+1}\right)^{-2}\right),~h\sim N(0,1),\ \textrm{for} \ l=7, \dots, 12.
\end{eqnarray*}
%\citet{collazos2016consistent}  considered  only the first six functions as the predictors. For simplicity, we convert each interval of definition, $\mathbb{T}_l$, into $[0,~1]$. To introduce the correlation between each functional variables,
%the random data $G_{il}$ were converted in to the function data $x_{il}$ through the linear transformations as follows:
%\begin{eqnarray*}
%	x_{i1}(t)&=&\sqrt{0.84} Z_{i1}(t)+0.4  Z_{i6}(t),\ \ \  	x_{i2}(t)=\sqrt{0.98} Z_{i2}(t)+0.1 Z_{i1}(t)+0.1  Z_{i5}(t),\\
%	x_{i3}(t)&=&\sqrt{0.84} Z_{i3}(t)+0.4  Z_{i4}(t),\ \ \   x_{i5}(t)=\sqrt{0.99} Z_{i5}(t)+0.1 Z_{i2}(t),\\
%	x_{il}(t)&=&Z_{il}(t)\ \ \ \ \ \textrm{for}\ \ l=4, 6, 7,\dots, 12.
%\end{eqnarray*}
The functional coefficients $\bm \beta(t)$ are generated based on the following 4 functions:
\begin{eqnarray*}
	f_1(t)&=& .03f(t,20,60)-.05f(t,50,20),\\
	f_2(t)&=& 4\sin(4\pi x)-\textrm{sign}(x-.3)-\textrm{sign}(.72-x),\\
	f_3(t)&=& -3\cos (2\pi t)+3 e^{t^2} / (t^3+1),\\
	f_4(t)&=&.1\sin(2\pi t)+.2\cos(2\pi t)+.3 \sin^2(2\pi t)+.4 \cos^3(2\pi t)+.5 \sin^3(2\pi t),
\end{eqnarray*}
where $f(t,\alpha,\beta)$ is the density function for beta distribution: $\rm{Beta}(\alpha,\beta)$.
%=\frac{\Gamma(\alpha+\beta)}{\Gamma(\alpha)\Gamma(\beta)}t^{\alpha-1}(1-t)^{\beta-1}.$
Note $f_1(t)$ has also been considered by \cite{zhao2012wavelet}; the second function $f_2$, the so-called ``Heavi-Sine'' function, is one of test functions from \cite{donoho1994ideal} which is very popular among wavelet literature (\citealp{antoniadis2001wavelet}); and $f_4$ was proposed by \citet{lin2013variable}.

\begin{figure}[htbp]
	\begin{center}
		\includegraphics[scale=.42]{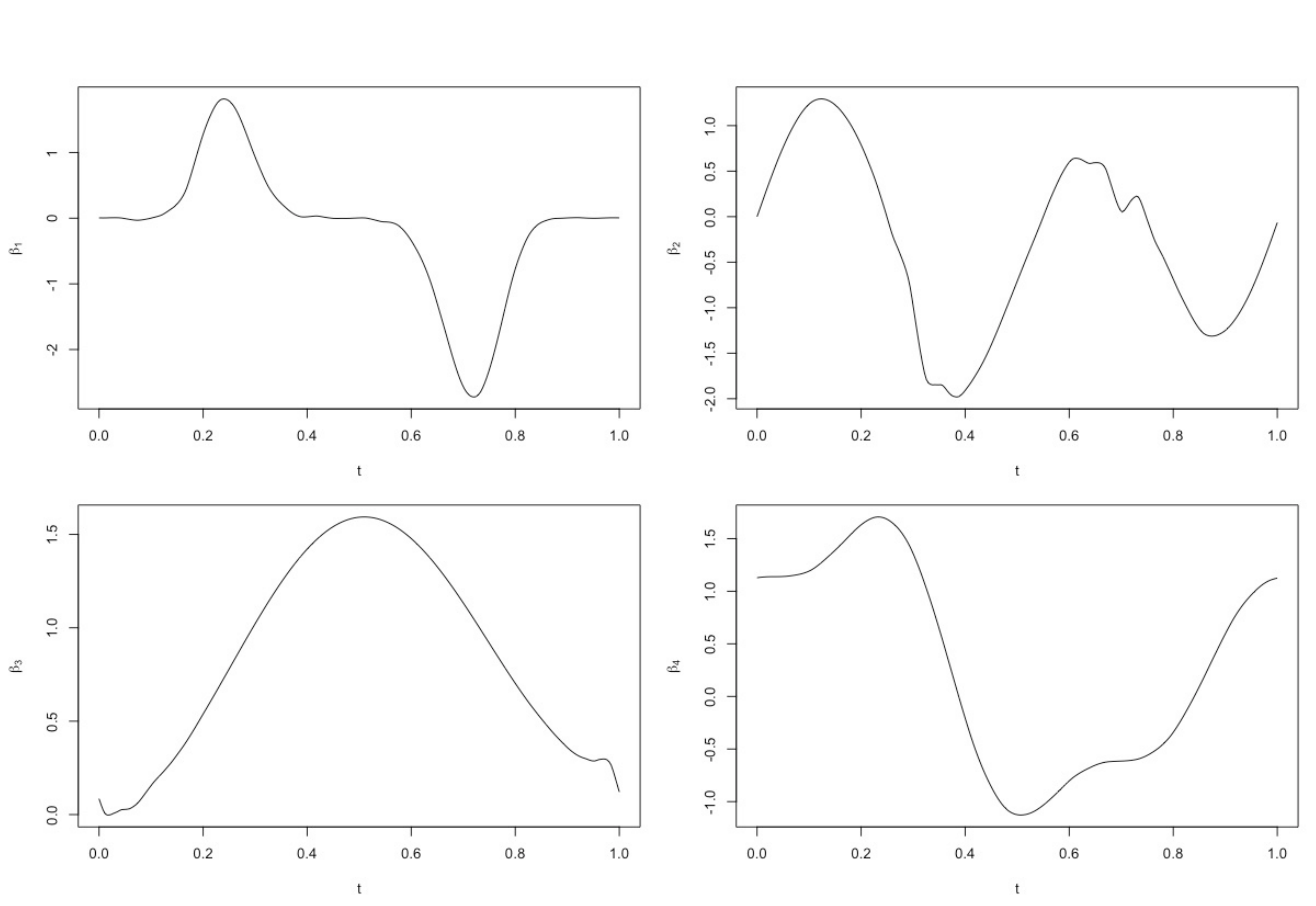}
		\caption{Slope functions of $\bm \beta_1$ to $\bm \beta_4$.}
		\label{slop_figure}
	\end{center}
\end{figure}

To generate the functional slopes $\beta_1(t),\dots \beta_4(t)$, we first apply DWT for $f_1,\dots, f_4$ and select the wavelet coefficients with absolute values greater than .1;
and based on the inverse DWT of the selected coefficients, we generate normalized $\beta_1(t),\dots \beta_4(t)$, each of which possesses sparsity and is shown in Figure \ref{slop_figure}.
The rest of slopes are set to be zero, i.e., $\beta_l(t)=0$ for $l=5,\ \dots,\ 12.$
The error term $\varepsilon_i$ is drawn from the following distributions:
%\begin{enumerate}
	1) Standard normal : $N(0,1)$;
	2) Mixed-variance:  $.95N(0,1)+.05N(0,10)$;
	3) t distribution with 3 degrees of freedom: $t_3$;
	4) Standard Cauchy: $C(0,1)$.
%\end{enumerate}
The signal-to-noise (SNR) ratio, defined as ${\mu}/{\sigma}$ in this paper, is chosen from three different levels: $\rm{SNR}=1,5,10$, where $\mu$ is the mean of signal and $\sigma$ is the standard deviation of the noise.
%Furthermore, we use signal-to-noise (SNR) ratio to measure the relativity of a desired signal to the background noise.
%While there are various definitions of it, in this paper, we choose to use $\rm{SNR}={\mu}/{\sigma}$, where $\mu$ is the mean of signal, and $\sigma$ is the standard deviation of the noise. The SNR is then set in three different levels such as $\rm{SNR}=1,5,10$.

The sizes of the training, tuning and testing data sets are $n$, $n$ and $10n$ respectively.
We select the tuning parameters via a grid search using the GIC and validation set methods through the tuning data set.
In GIC, $\phi_n$s are $5p_n$, $5p_n$ and $p_n$ for the quantile sparse group lasso, lasso and group lasso methods respectively, while  $p_n=\log\left(\log(n)\right)\log\left(\log(p)\right)/\left(10n\right)$. %In real data example later, we use exactly the same $\phi_n$ for these three methods.
The validation set method is used to select the gold standard (GS) tuning parameters that minimize the prediction error of tuning data sets (\citealp{li2007quantile}, \citealp{zou2008composite}, \citealp{wu2009variable}).
%Similarly defined testing errors on the testing data set are reported. More explicitly, a test error refers to the average check loss on the independent testing data set. Using these $\lambda_1,\ \lambda_2$'s, we calculated the prediction errors and the mean absolute deviation with test data for each criterion.

In our simulations, we choose $n=200$, $400$, set $\tau=0.5$, and use $100$ Monte Carlo repititions.
We use the following five criteria of the performance, namely,
the group accuracy (GA), variable accuracy (VA), mean absolute prediction error (MAPE), mean integrated square errors (MISE) and individual  integrated square errors (ISE). % of the estimated functional coefficients.
% the proportion of correctly picked up and dropped functional components, called
The group accuracy (GA) is the proportion of correctly picked up and dropped off functional components, that is $\textrm{GA}=E\left(\left({|\widehat{M}\cap M_0|+|\widehat{M}^c\cap M_0^c|}\right)/{12}\right)$ with
%$M_0$ and $\widehat{M}$ denote the set of indices of the true functional variables and selected functional variables, respectively, such as
$M_0=\{l:{\beta}_l(t)\ne 0\}$ and $\widehat{M}=\{l:\hat{\beta}_l(t)\ne 0\}$.
The variable accuracy (VA) is defined similarly as GA by simply replacing the $M_0$ and $\widehat{M}$ as the true and estimated index sets of non-zero wavelet coefficients. The mean absolute prediction error (MAPE) is MAPE$=E(|\hat{y}-y|)$. % is assessed on the testing data set for each Monte Carlo repetition.
The mean integrated square errors (MISE) of the $12$ estimated functional coefficients:
\begin{equation*}
	\textrm{MISE} = {\frac 1 {12} \sum_{l=1}^{12} \int_0^1 (\hat{\beta}_l(t)-\beta_l(t))^2}dt,
\end{equation*}
%\begin{equation*}
%	\textrm{MISE} = {\frac 1 {12} \sum_{l=1}^{12}\sum_{j=1}^{256}(\hat{\beta}_l(t_j)-\beta_l(t_j))^2},
%\end{equation*}
as well as the individual integrated square error (ISE):
\begin{equation*}
	\textrm{ISE}_l = {\int_0^1 (\hat{\beta}_l(t)-\beta_l(t))^2}dt,
\end{equation*}
is used to measure the estimation accuracy of functional coefficients.

Due to space limit, we only discuss the results of SNR $=5$. The results for the other two SNRs are both in favor of our method and deferred to the Appendix.
%The simulation results for SNR =  5 are summarized in Table \ref{tab:error:snr5} and \ref{tab:function error:snr5}. %As shown in Table \ref{tab:error:snr5} and \ref{tab:function error:snr5},
As shown in Table \ref{tab:error:snr5}, in general, the performance of qSGL method is better than the qL and qGL methods in terms of mean integrated square errors (MISEs) and mean absolute prediction errors (MAPEs).
For different error types, our proposed GIC approach is only slightly outperformed by the gold standards.
As the sample size increases, the MISEs and MAPEs decrease, which is consistent with our theoretical results.
For group accuracy (GA), qGL performs better than the other methods in most cases, while qL performs quite well in terms of variable accuracy (VA).
However, in the case of GIC,  the sparse group lasso method outperforms the two competitors regarding both GA and VA, especially for larger sample sizes.
In Table \ref{tab:function error:snr5}, it shows that the ISEs of sparse group lasso are smaller than the other two methods.
It also shows that the ISE of $\hat{\beta}_1(t)$ is always less than the other three slope functions in most cases regardless the methods used.
It might be due to the fact that $\beta_1(t)$ is smoother than the other slopes; see Figure \ref{slop_figure}.
\begin{table}[htbp]
	\centering
	\small
	\begin{tabular}{|r|r|r|rrrr|rrrr|}
		\hline
		% Sample size& Noise Type	&Method	& MISE & Group Accuracy  & Variable Accuracy & Prediction Error & MISE & Group Accuracy & Variable Accuracy & Prediction Error \\
		&  & 	  &  &GS  &  &  &  & GIC &  &  \\
		n& Noise&Method	& MISE & GA  & VA & MAPE & MISE & GA & VA & MAPE \\
		\hline
		&    &qSGL & 1.449 & 0.930 & 0.934 & 2.600 & 1.522 & 0.594 & 0.840 & 2.851 \\
		&1  &qL & 3.230 & 0.919 & 0.961 & 2.871 & 3.159 & 0.482 & 0.904 &3.205 \\
		& 	&qGL & 1.835 & 1.000 & 0.082 & 2.862 & 2.121 & 0.970 & 0.343 & 4.763 \\
		\cline{2-11}
		& 	&qSGL & 1.372 & 0.960 & 0.934 & 2.466 & 1.516 & 0.623 & 0.835 & 2.796 \\
		&2 	&qL & 3.023 & 0.932 & 0.960 & 2.749 & 3.086 & 0.496 & 0.905 & 3.142 \\
		& 	&qGL & 1.802 & 1.000 & 0.082 & 2.781 & 2.068 & 0.973 & 0.326 & 4.476 \\
		\cline{2-11}
		200 & 	&qSGL & 0.598 & 1.000 & 0.911 & 1.436 & 0.932 & 0.871 & 0.857 & 1.953 \\
		&3	&qL & 1.420 & 0.985 & 0.945 & 1.671 & 2.487 & 0.686 & 0.909 & 2.654 \\
		& 	&qGL & 1.630 & 1.000 & 0.065 & 2.386 & 1.735 & 0.993 & 0.140 & 2.836 \\
		\cline{2-11}
		& 	&qSGL & 1.284 & 0.972 & 0.934 & 2.326 & 1.497 & 0.617 & 0.829 & 2.755 \\
		&4   &qL & 2.826 & 0.927 & 0.958 & 2.625 & 3.135 & 0.490 & 0.907 & 3.145 \\
		& 	&qGL & 1.775 & 1.000 & 0.075 & 2.656 & 2.043 & 0.976 & 0.295 & 4.225 \\
		\hline
		& 	&qSGL & 0.925 & 0.989 & 0.915 & 2.095 & 1.224 & 0.911 & 0.920 & 2.220 \\
		&1 &qL & 1.774 & 0.944 & 0.946 & 2.187 & 2.125 & 0.617 & 0.898 & 2.371 \\
		& 	&qGL & 1.581 & 1.000 & 0.054 & 2.393 & 2.246 & 0.958 & 0.569 & 5.240 \\
		\cline{2-11}
		& 	&qSGL & 0.842 & 0.995 & 0.911 & 1.954 & 1.105 & 0.967 & 0.937 & 2.058 \\
		&2&qL & 1.640 & 0.965 & 0.947 & 2.040 & 1.853 & 0.729 & 0.912 & 2.190 \\
		& 	&qGL & 1.549 & 1.000 & 0.056 & 2.306 & 2.263 & 0.957 & 0.582 & 5.294 \\
		\cline{2-11}
		400&	&qSGL & 0.157 & 1.000 & 0.875 & 1.001 & 0.272 & 1.000 & 0.930 & 1.108 \\
		
		&3&qL & 0.285 & 1.000 & 0.908 & 1.026 & 0.481 & 0.991 & 0.943 & 1.108 \\
		& 		&qGL & 1.255 & 1.000 & 0.050 & 1.996 & 1.438 & 0.992 & 0.155 & 2.472 \\
		\cline{2-11}
		& 	&qSGL & 0.738 & 0.996 & 0.909 & 1.785 & 0.995 & 0.983 & 0.939 & 1.910 \\
		&4&qL & 1.469 & 0.978 & 0.947 & 1.860 & 1.737 & 0.735 & 0.906 & 2.052 \\
		& 	&qGL & 1.505 & 0.999 & 0.054 & 2.194 & 2.102 & 0.969 & 0.499 & 4.490 \\
		\hline
	\end{tabular}
	\caption{Simulation summary of SNR=5. The first column $n$ is the size of training data. The second column is the type of noise. The third column is the method we used, qSGL for the quantile sparse group lasso, qL for the quantile Lasso, and qGL for the quantile group lasso. GS means $\lambda$ was selected by the validation method (gold standard). GIC means $\lambda$ selected via the GIC criterion. MISE stands for mean integrated errors. MAPE, GA and VA indicate mean absolute prediction error, group accuracy and variable accuracy, respectively.}
	\label{tab:error:snr5}
\end{table}

\begin{table}[htbp]
	\centering
	\small
	\begin{tabular}{|r|r|r|rrrr|rrrr|}
		\hline
		&  & 	  &  &GS  &  &  &  & GIC &  &  \\
		n& Noise&Method	& ISE1  & ISE2  & ISE3 & ISE4 & ISE1 & ISE2 & ISE3 & ISE4 \\
		\hline
		&    &qSGL & 0.116 & 0.585 & 0.331 & 0.385 & 0.133 & 0.550 & 0.322 & 0.387 \\
		&1  &qL & 0.289 & 0.758 & 1.386 & 0.734 & 0.318 & 0.618 & 1.136 & 0.732 \\
		& 	&G & 0.351 & 0.675 & 0.359 & 0.447 & 0.372 & 0.728 & 0.370 & 0.648 \\
		\cline{2-11}
		& 	&qSGL & 0.116 & 0.540 & 0.322 & 0.368 & 0.137 & 0.560 & 0.318 & 0.377 \\
		&2 	&qL & 0.283 & 0.674 & 1.302 & 0.703 & 0.336 & 0.631 & 1.049 & 0.740 \\
		& 	&qGL & 0.348 & 0.665 & 0.349 & 0.438 & 0.367 & 0.714 & 0.362 & 0.621 \\
		\cline{2-11}
		200 & 	&qSGL & 0.051 & 0.162 & 0.163 & 0.214 & 0.077 & 0.311 & 0.221 & 0.267 \\
		&3	&qL & 0.105 & 0.204 & 0.614 & 0.468 & 0.238 & 0.460 & 0.939 & 0.610 \\
		& 	&qGL &  0.332 & 0.605 & 0.297 & 0.395 & 0.342 & 0.632 & 0.313 & 0.446 \\
		\cline{2-11}
		& 	&qSGL & 0.104 & 0.498 & 0.304 & 0.354 & 0.129 & 0.551 & 0.328 & 0.367 \\
		&4   &qL & 0.248 & 0.628 & 1.211 & 0.679 & 0.318 & 0.613 & 1.157 & 0.707 \\
		& 	&qGL & 0.345 & 0.657 & 0.343 & 0.427 & 0.367 & 0.709 & 0.366 & 0.597 \\
		\hline
		& 	&qSGL & 0.074 & 0.321 & 0.217 & 0.293 & 0.091 & 0.470 & 0.265 & 0.353 \\
		&1 &qL & 0.141 & 0.318 & 0.729 & 0.532 & 0.155 & 0.377 & 0.719 & 0.575 \\
		& 	&qGL &  0.325 & 0.590 & 0.285 & 0.381 & 0.363 & 0.731 & 0.399 & 0.752 \\
		\cline{2-11}
		& 	&qSGL & 0.071 & 0.274 & 0.207 & 0.273 & 0.088 & 0.421 & 0.248 & 0.331 \\
		&2&qL &  0.117 & 0.279 & 0.695 & 0.508 & 0.139 & 0.324 & 0.675 & 0.519 \\
		& 	&qGL & 0.321 & 0.577 & 0.278 & 0.372 & 0.364 & 0.736 & 0.401 & 0.761 \\
		\cline{2-11}
		400&	&qSGL & 0.010 & 0.018 & 0.063 & 0.065 & 0.016 & 0.045 & 0.094 & 0.115 \\
		&3&qL & 0.012 & 0.017 & 0.139 & 0.110 & 0.018 & 0.034 & 0.234 & 0.187 \\
		& 		&G &0.295 & 0.446 & 0.205 & 0.308 & 0.311 & 0.504 & 0.244 & 0.375 \\
		\cline{2-11}
		& 	&qSGL & 0.057 & 0.220 & 0.195 & 0.253 & 0.071 & 0.366 & 0.233 & 0.312 \\
		&4&qL & 0.096 & 0.218 & 0.643 & 0.478 & 0.116 & 0.273 & 0.631 & 0.515 \\
		& 	&qGL & 0.319 & 0.555 & 0.266 & 0.363 & 0.354 & 0.700 & 0.383 & 0.664 \\
		\hline
	\end{tabular}
	\caption{Individual funtional $L_2$ error of SNR=5. The first column n is the size of training data. The second column is the noise type. The third column is the method we used. ISE1: $||\hat{\beta}_1-\beta_1||_2^2$; ISE2: $||\hat{\beta}_2-\beta_2||_2^2$; ISE3: $||\hat{\beta}_3-\beta_3||_2^2$; ISE4: $||\hat{\beta}_4-\beta_4||_2^2$. }
	\label{tab:function error:snr5}
\end{table}

\subsection{Real Data}
The real data we use is a subset of the ADHD-200 Sample Initiative Project (\citealp{mennes2013making}), which studies attention deficit hyperactivity disorder (ADHD), the most commonly diagnosed mental disorder of childhood which may persist into adulthood.
ADHD is characterized by problems related to paying attention, hyperactivity, or impulsive behavior.
The dataset is a filtered preprocessed resting state fMRI data from New York University Child Study Centre using the Anatomical Automatic Labeling (AAL) atlas (\citealp{tzourio2002automated}). In the dataset, there are 172 equally spaced time courses in the filtering and AAL contains 116 Regions of Interests (ROIs) fractionated into functional space using nearest-neighbor interpolation. Each of 172 time courses is then smoothed to 64 equally to apply DWT.
After cleaning the raw data that fails in quality control or has missing data, we have 120 individuals in final analysis.
Grouping ROIs in terms of their anatomical functions and averaging within each group the corresponding time courses, we have 59 averaged time courses of grouped ROIs serving as functional predictors, each of which has 64 equally spaced time points.
In addition, 8 scalar covariates are considered, including gender, age, handedness, diagnosis status, medication status, Verbal IQ, Performance IQ and Full4 IQ. The response of interest is the ADHD index, a measurement of severity of mental disorder.

We apply partial functional linear quantile regression model  (\ref{quantile:model}) with $59$ functional covariates and $8$ scalar covariates.
In order to select the significant functional covariates from $59$ ROIs,
we use the procedure proposed by \cite{Meinshausen2010} to obtain stable selections from $100$ bootstrap samples.
The tuning parameters are chosen by GIC.
%proposed functional variable selection method and compare the results with those of the quantile Lasso and the quantile group lasso selection procedures by using the same wavelet basis functions. First,
%The $L_2$ norms of the estimated slope functions from bootstrap samples are calculated.
The boxplots of $L_2$ norms of the estimated slope functions from bootstrap samples are shown in Figure \ref{fig1:boxplot:sgl}, \ref{fig1:boxplot:l} and \ref{fig1:boxplot:g} in the Appendix.
The selection criterion is that the median of corresponding $L_2$ norm should be greater than $ 10^{-5}$.
%Finally, the empirical distributions of each slope function's $L_2$ norm are estimated from this bootstrapping. The boxplots of the $L_2$ norms for each method are shown in the Appendix; see Figure \ref{fig1:boxplot:sgl}, \ref{fig1:boxplot:l}, \ref{fig1:boxplot:g}. After sorting medians of each function's $L_2$ norm, we select the important functional predictors by the threshold of $10^{-5}$.

In neurological science literature on ADHD, it has been shown that the $7$ regions of cerebellum, temporal, vermis, parietal, occipital, cingulum and frontal are commonly discovered to be significantly related to ADHD symptoms  from various studies ( \citealp{max2005prefrontal}; \citealp{konrad2010adhd}; \citealp{tomasi2012abnormal}).
We first evaluate the performances of qSGL, qL and qGL methods in terms of the selection of these 7 regions, which are essentially 14 ROIs including the left and right parts.
In Table \ref{Result:adhd} and \ref{Result:adhd1}, we list the selected ROIs from three different methods. In particular, qSGL, qL and qGL select 15, 20 and 9 ROIs respectively.
In terms of those 7/14 commonly discovered regions/ROIs,
Both our proposed qSGL and qGL methods have lower false discovery rates ($33\%$) than the qL method ($55\%$),
while our method is superior to the qGL as it identifies more true positives ($10$ vs $6$).
Moreover, ``Occipital R'', the right occipital region, can only be identified by our method.
While both Table  \ref{Result:adhd} and \ref{Result:adhd1} confirm that most of the selected ROIs are coming from the 7/14 mostly discovered regions/ROIs,  the three methods also suggest three other common ROIs: ``Olfactory R'', ``Supramarginal R'', and ``Caudate R'', namely
right olfactory, right supramarginal, and right caudate regions respectively,
which have been evidently important as suggested by some ADHD studies. %  There has been some medical evidence to suggest that the affected nodes include these other three regions.
For instance, \citet{schrimsher2002caudate} revealed a relationship between caudate asymmetry and some symptoms related to ADHD.
%This correlation is congruent with previous associations of the caudate with attentional functioning.
The findings of \cite{Sidlauskaite2015506} imply the supramarginal gyrus is associated with the ADHD symptom scores.

\begin{table}[htbp]
	\centering
	\small
	\begin{tabular}{|c|c|}
		\hline
		% Sample size& Noise Type	&Method	& MISE & Group Accuracy  & Variable Accuracy & Prediction Error & MISE & Group Accuracy & Variable Accuracy & Prediction Error \\
		Method&    Significant~ROIs \\
		
		\hline
		
		&``Temporal R'' ~``Cerebelum R" ~  ``Frontal R"  ~ ``Occipital R" ~ ``Olfactory R" \\

		qSGL & ``SupraMarginal R"  ~ ``Caudate R"       ~  ``Vermis"            ~``Cuneus L"   ~   ``Parietal R"     \\
		&  ``Frontal L"     ~    ``Precuneus R"
		~ ``Temporal L"   ~    ``Cerebelum L"      ~ ``Precentral R" \\
		
		\hline	
		& ``Frontal R"   ~      ``Caudate R"    ~     ``Temporal R"
		~ ``Cuneus L"      ~   ``SupraMarginal R"   \\
		& ``Parietal R"       ~  ``Lingual L"      ~   ``Frontal L"         ~ ``Precuneus R" ~
		``Vermis"    \\
		qL	& ``Fusiform R"     ~   ``Pallidum L"
		~ ``Olfactory R"     ~  ``Precentral R"    ~  ``Cingulum L"  \\
		
		& ``Cuneus R"      ~    ``Parietal L"      ~  ``Temporal L"
		~ ``Angular L"       ~ ``Cerebelum R" \\

		\hline

		& ``Caudate R" ~     ``Frontal R"     ~    ``Cerebelum R"   ~
		``Vermis"      ~ ``Olfactory R"   \\
		qGL 	    & ``Temporal R"  ~    ``Precentral R"    ~ ``SupraMarginal R"  ~ ``Frontal L"     \\
		
		\hline
		
	\end{tabular}
	\caption{Selected ROIs for the  ADHD-200 fMRI Dataset.}
	\label{Result:adhd}
\end{table}

\begin{table}[htbp]
	\centering
	\begin{tabular}{|c|c|c|c|}
		\hline
		Significant~regions&  qSGL & qL & qGL  \\		
		\hline		
		Cerebellum &   R L   &   R     &    R    \\
		\hline
		Temporal &    R L  &  R L     &    R    \\
		\hline		
		Vermis&    R L  &   R L     &   R L     \\
		\hline	
		Parietal & R     &  R L      &        \\
		\hline
		Occipital&   R   &      &       \\
		\hline
		Cingulum&  &    L    &        \\
		\hline
		Frontal&   R L   &   L     &    R L    \\
		\hline
	\end{tabular}
	\caption{Selected ROIs for the suggested 7 regions, `R'  and `L' indicate the region is selected from the right brain and left brain, respectively.  Blank means the brain region is not chosen.}
	\label{Result:adhd1}
\end{table}

%Selected predictors for the suggestion 7 regions, R and L indicates, the region has been selected from the right brain or left brain, respectively. The symbol $\times$ means the brain region has not been chosen,

\section{Discussion}
This article studies quantile regression in partial functional linear model where response is scalar and predictors include both scalars and multiple functions.
We adopt wavelet basis to well approximate functional slopes while effectively detect local features.
A sparse group lasso method is proposed to select important functional predictors while capture shared information among them.
We reformulate the proposed problem into a standard second-order cone program and then solve it by an interior point method.
A novel and efficient algorithm by using alternating direction method of multipliers (ADMM) is utilized
to solve the optimization problem.
In addition, we successfully derive the asymptotic properties including the convergence rate and prediction error bound which guarantee a good theoretical performance of the proposed method.
Simulation studies demonstrate that our proposed method is more effective in estimating coefficients and making predictions while capable of identifying non-zero functional components and wavelet coefficients.
We analyze a real data from ADHD-200 fMRI data set and show the superiority of our method.  Moreover, our analysis makes some new discovery about other brain regions that are evidently important in making diagnosis.

There are several topics that merit further research. Other asymptotic properties, such as the model selection consistency and asymptotic normality, of our proposed method could be developed. The technique proposed to reformulate our problem into a second order cone program (SOCP) could be further adapted to other penalized quantile regression problems; for example, quantile ridge regression (\citealp{wu2009variable}). Moreover, to estimate the functional slopes, the wavelet-based technique can also be used together with principal component analysis or partial least squares methods (\citealp{reiss2015wavelet}).

\section{Appendix}

\begin{table}[htbp]
	\centering
	\small
	\begin{tabular}{|r|r|r|rrrr|rrrr|}
		\hline
		% Sample size& Noise Type	&Method	& MISE & Group Accuracy  & Variable Accuracy & Prediction Error & MISE & Group Accuracy & Variable Accuracy & Prediction Error \\
		&  & 	  &  &GS  &  &  &  & GIC &  &  \\
		n& Noise&Method	& MISE & GA  & VA & MAPE & MISE & GA & VA & MAPE \\
		\hline
		&    &qSGL & 2.426 & 0.860 & 0.959 & 9.557 & 6.361 & 0.480 & 0.854 & 11.720 \\
		& 1   &qL & 5.885 & 0.965 & 0.972 & 9.553 & 17.062 & 0.358 & 0.891 & 13.134 \\
		& 	&qGL & 2.601 & 0.852 & 0.118 & 9.637 & 4.091 & 0.708 & 0.406 & 12.728 \\
		\cline{2-11}
		& 	&qSGL & 2.322 & 0.876 & 0.958 & 8.833 & 6.013 & 0.509 & 0.857 & 10.968 \\
		&2 	&qL & 5.592 & 0.968 & 0.971 & 8.844 & 16.564 & 0.363 & 0.891 & 12.760 \\
		& 	&qGL & 2.619 & 0.870 & 0.123 & 8.973 & 4.473 & 0.704 & 0.374 & 11.844 \\
		\cline{2-11}

		200 & 	&qSGL & 1.063 & 0.994 & 0.930 & 4.200 & 1.594 & 0.891 & 0.908 & 4.774 \\
		&3	&qL & 2.462 & 0.978 & 0.958 & 4.491 & 7.252 & 0.547 & 0.911 & 7.466 \\
		& 	&qGL & 1.741 & 1.000 & 0.073 & 4.776 & 3.699 & 0.857 & 0.330 & 7.875 \\
		\cline{2-11}
		
		& 	&qSGL & 2.252 & 0.925 & 0.958 & 7.967 & 5.795 & 0.510 & 0.856 & 10.353 \\
		&4   &qL& 5.332 & 0.983 & 0.971 & 8.012 & 15.874 & 0.365 & 0.891 & 12.401 \\
		& 	&qGL & 2.402 & 0.920 & 0.113 & 8.099 & 4.152 & 0.751 & 0.404 & 11.165 \\
		\hline

		& 	&qSGL & 2.186 & 0.935 & 0.954 & 8.699 & 2.427 & 0.959 & 0.974 & 9.529 \\
		&1 &qL & 5.246 & 0.981 & 0.971 & 8.756 & 5.916 & 0.966 & 0.970 & 8.906 \\
		& 	&qGL & 2.336 & 0.944 & 0.106 & 8.788 & 3.450 & 0.877 & 0.667 & 11.703 \\
		
		\cline{2-11}

		& 	&qSGL & 2.126 & 0.954 & 0.954 & 8.083 & 2.414 & 0.963 & 0.976 & 9.030 \\
		&2&qL & 4.962 & 0.983 & 0.970 & 8.153 & 5.175 & 1.000 & 0.974 & 8.206 \\
		& 	&qGL & 2.234 & 0.973 & 0.102 & 8.182 & 2.742 & 0.898 & 0.718 & 11.403 \\
		
		\cline{2-11}

		400&	&qSGL &0.492 & 1.000 & 0.883 & 3.630 & 1.004 & 0.999 & 0.951 & 3.985 \\
		&3&qL & 1.035 & 0.995 & 0.934 & 3.698 & 1.855 & 0.994 & 0.965 & 4.018 \\
		& 		&qGL & 1.415 & 1.000 & 0.052 & 4.305 & 2.394 & 0.932 & 0.551 & 7.679 \\
		
		\cline{2-11}

		& 	&qSGL & 2.008 & 0.962 & 0.950 & 7.301 & 2.338 & 0.965 & 0.975 & 8.258 \\
		&4&qL & 4.602 & 0.983 & 0.970 & 7.394 & 5.991 & 0.967 & 0.968 & 7.634 \\
		& 	&qGL & 2.133 & 0.983 & 0.102 & 7.376 & 3.250 & 0.888 & 0.692 & 10.880 \\
		
		\hline
	\end{tabular}
	%\caption*{ As for Table \ref{tab:error:snr5}, 	
	\caption{Simulation summary of SNR=1, as for Table \ref{tab:error:snr5}.}
	\label{tab:error:snr1}
\end{table}

\begin{table}[htbp]
	\centering
	\small
	\begin{tabular}{|r|r|r|rrrr|rrrr|}
		\hline
		&  & 	  &  &GS  &  &  &  & GIC &  &  \\
		n& Noise&Method	& ISE1  & ISE2  & ISE3 & ISE4 & ISE1 & ISE2 & ISE3 & ISE4 \\
		\hline
		&    &qSGL  & 0.186 & 0.822 & 0.673 & 0.684 & 0.270 & 2.415 & 0.749 & 0.654 \\
		&1  &qL &  0.629 & 1.004 & 3.197 & 0.987 & 1.073 & 3.883 & 3.019 & 2.129 \\
		& 	&qGL & 0.407 & 0.901 & 0.537 & 0.693 & 0.529 & 1.851 & 0.555 & 0.868 \\
		\cline{2-11}
		
		& 	&qSGL & 0.181 & 0.810 & 0.642 & 0.635 & 0.264 & 2.494 & 0.712 & 0.640 \\
		&2 	&qL & 0.592 & 0.973 & 2.974 & 0.989 & 1.078 & 4.370 & 2.808 & 1.825 \\
		& 	&qGL&0.411 & 0.971 & 0.521 & 0.660 & 0.530 & 2.198 & 0.583 & 0.861 \\
		\cline{2-11}
		
		200 & 	&qSGL & 0.087 & 0.394 & 0.252 & 0.315 & 0.112 & 0.646 & 0.322 & 0.383 \\
		&3	&qL & 0.197 & 0.520 & 1.067 & 0.640 & 0.575 & 1.677 & 1.619 & 1.045 \\
		& 	&qGL & 0.342 & 0.645 & 0.330 & 0.422 & 0.438 & 1.729 & 0.497 & 0.737 \\
		\cline{2-11}

		& 	&qSGL & 0.165 & 0.816 & 0.646 & 0.589 & 0.243 & 2.383 & 0.764 & 0.641 \\
		&4   &qL& 0.552 & 0.961 & 2.781 & 0.986 & 0.982 & 4.010 & 2.891 & 1.769 \\
		& 	&qGL & 0.396 & 0.858 & 0.511 & 0.605 & 0.509 & 1.989 & 0.544 & 0.862 \\
		\hline
		
		& 	&qSGL & 0.163 & 0.830 & 0.593 & 0.565 & 0.189 & 0.801 & 0.616 & 0.817 \\
		&1 &qL & 0.565 & 0.973 & 2.692 & 0.966 & 0.549 & 1.176 & 2.773 & 1.060 \\
		& 	&qGL &  0.387 & 0.837 & 0.492 & 0.598 & 0.453 & 1.339 & 0.422 & 1.041 \\
		
		\cline{2-11}

		& 	&qSGL & 0.165 & 0.814 & 0.579 & 0.540 & 0.194 & 0.795 & 0.619 & 0.803 \\
		&2&qL & 0.513 & 0.966 & 2.501 & 0.938 & 0.523 & 0.970 & 2.679 & 0.982 \\
		& 	&qGL & 0.383 & 0.797 & 0.476 & 0.563 & 0.420 & 0.900 & 0.375 & 1.007 \\
		\cline{2-11}

		400& &qSGL & 0.038 & 0.133 & 0.137 & 0.177 & 0.070 & 0.393 & 0.241 & 0.298 \\
		&3&qL & 0.065 & 0.147 & 0.456 & 0.350 & 0.123 & 0.404 & 0.814 & 0.500 \\
		&  &qGL & 0.312 & 0.516 & 0.242 & 0.344 & 0.383 & 0.778 & 0.405 & 0.802 \\
		
		\cline{2-11}
		
		& 	&qSGL &0.146 & 0.794 & 0.540 & 0.502 & 0.176 & 0.799 & 0.604 & 0.758 \\
		&4&qL &  0.414 & 0.952 & 2.281 & 0.919 & 0.461 & 1.269 & 2.478 & 1.140 \\
		& 	&qGL & 0.376 & 0.771 & 0.448 & 0.527 & 0.433 & 1.177 & 0.412 & 1.041 \\
		\hline
	\end{tabular}
	\caption{Individual funtional $L_2$ error when SNR=1, as for Table \ref{tab:function error:snr5}.}
	\label{tab:function error:snr1}
\end{table}

\begin{table}[htbp]
	\centering
	\small
	\begin{tabular}{|r|r|r|rrrr|rrrr|}
		\hline
		% Sample size& Noise Type	&Method	& MISE & Group Accuracy  & Variable Accuracy & Prediction Error & MISE & Group Accuracy & Variable Accuracy & Prediction Error \\
		&  & 	  &  &GS  &  &  &  & GIC &  &  \\
		n& Noise&Method	& MISE & GA  & VA & MAPE & MISE & GA & VA & MAPE \\
		\hline
		
		&    &qSGL & 0.907 & 0.988 & 0.906 & 1.617 & 0.920 & 0.935 & 0.839 & 1.683 \\
		& 1   &qL & 1.962 & 0.917 & 0.939 & 1.835 & 1.964 & 0.792 & 0.910 & 1.917 \\
		& 	&qGL &  1.679 & 1.000 & 0.064 & 2.195 & 1.743 & 0.994 & 0.132 & 2.578 \\
		
		\cline{2-11}
		& 	&qSGL & 0.898 & 0.992 & 0.912 & 1.576 & 0.913 & 0.943 & 0.840 & 1.662 \\
		&2 	&qL & 1.866 & 0.932 & 0.942 & 1.784 & 1.917 & 0.790 & 0.912 & 1.888 \\
		& 	&qGL & 1.669 & 1.000 & 0.067 & 2.172 & 1.779 & 0.989 & 0.161 & 2.857 \\
		\cline{2-11}
		
		200 & 	&qSGL & 0.498 & 1.000 & 0.903 & 1.124 & 0.709 & 0.943 & 0.849 & 1.482 \\
		&3	&qL & 1.203 & 0.993 & 0.943 & 1.325 & 1.756 & 0.828 & 0.914 & 1.867 \\
		& 	&qGL & 1.603 & 1.000 & 0.062 & 2.170 & 1.659 & 0.995 & 0.109 & 2.465 \\
		\cline{2-11}
		& 	&qSGL & 0.842 & 0.992 & 0.915 & 1.502 & 0.911 & 0.943 & 0.843 & 1.656 \\
		&4   &qL& 1.774 & 0.952 & 0.944 & 1.709 & 1.928 & 0.792 & 0.913 & 1.904 \\
		& 	&qGL & 1.656 & 1.000 & 0.065 & 2.116 & 1.722 & 0.996 & 0.125 & 2.420 \\
		\hline
		
		& 	&qSGL &  0.499 & 0.999 & 0.892 & 1.142 & 0.610 & 0.963 & 0.874 & 1.222 \\
		&1 &qL & 0.981 & 0.965 & 0.932 & 1.187 & 1.029 & 0.838 & 0.879 & 1.278 \\
		& 	&qGL &1.371 & 1.000 & 0.051 & 1.684 & 1.557 & 0.998 & 0.208 & 2.183 \\
		
		\cline{2-11}
		& 	&qSGL & 0.458 & 1.000 & 0.890 & 1.069 & 0.565 & 0.981 & 0.897 & 1.145 \\
		&2&qL & 0.902 & 0.975 & 0.933 & 1.114 & 0.927 & 0.867 & 0.894 & 1.190 \\
		& 	&qGL & 1.361 & 1.000 & 0.052 & 1.665 & 1.567 & 0.996 & 0.216 & 2.275 \\
		\cline{2-11}
		
		400&	&qSGL & 0.096 & 1.000 & 0.874 & 0.602 & 0.167 & 1.000 & 0.918 & 0.671 \\
		&3&qL & 0.151 & 1.000 & 0.903 & 0.617 & 0.299 & 0.999 & 0.941 & 0.681 \\
		& 		&qGL & 1.220 & 1.000 & 0.050 & 1.679 & 1.260 & 1.000 & 0.081 & 1.759 \\
		
		\cline{2-11}
		
		& 	&qSGL & 0.410 & 1.000 & 0.891 & 0.981 & 0.515 & 0.978 & 0.899 & 1.067 \\
		&4&qL & 0.837 & 0.988 & 0.934 & 1.025 & 0.866 & 0.898 & 0.898 & 1.105 \\
		& 	&qGL &  1.336 & 1.000 & 0.050 & 1.627 & 1.494 & 0.997 & 0.175 & 2.075 \\
		\hline
	\end{tabular}
	\caption{Simulation summary of SNR=10, as for Table \ref{tab:error:snr5}.}
	%\caption*{ As for Table \ref{tab:error:snr5}, but SNR=10}
	%{ The first column n is the size of training data. The second column is the type of Noise. The third column is the method we used; Q: quantile sparse group lassor; L: quantile lassor, G: quantile group lassor. GS means the $\lambda$ selected by the validation method, said the gold standard. GIC means the $\lambda$ selected by the GIC criterion,  }
	\label{tab:error:snr10}
\end{table}

\begin{table}[htbp]
	\centering
	\small
	\begin{tabular}{|r|r|r|rrrr|rrrr|}
		\hline
		&  & 	  &  &GS  &  &  &  & GIC &  &  \\
		n& Noise&Method	& ISE1  & ISE2  & ISE3 & ISE4 & ISE1 & ISE2 & ISE3 & ISE4 \\
		\hline
		&    &qSGL  & 0.080 & 0.298 & 0.220 & 0.286 & 0.082 & 0.292 & 0.222 & 0.284 \\
		&1  &qL & 0.166 & 0.340 & 0.819 & 0.570 & 0.165 & 0.317 & 0.799 & 0.569 \\
		& 	&qGL & 0.334 & 0.625 & 0.312 & 0.407 & 0.338 & 0.637 & 0.318 & 0.449 \\
		\cline{2-11}
		& 	&qSGL &0.081 & 0.299 & 0.216 & 0.282 & 0.087 & 0.294 & 0.218 & 0.277 \\
		&2 	&qL & 0.158 & 0.315 & 0.776 & 0.559 & 0.177 & 0.321 & 0.746 & 0.565 \\
		& 	&qGL& 0.334 & 0.621 & 0.310 & 0.403 & 0.342 & 0.641 & 0.318 & 0.477 \\
		\cline{2-11}
		
		200 & 	&qSGL & 0.040 & 0.117 & 0.146 & 0.188 & 0.061 & 0.206 & 0.182 & 0.233 \\
		&3	&qL & 0.077 & 0.148 & 0.540 & 0.415 & 0.141 & 0.265 & 0.737 & 0.512 \\
		& 	&qGL & 0.330 & 0.597 & 0.289 & 0.387 & 0.333 & 0.607 & 0.296 & 0.423 \\
		\cline{2-11}
		
		& 	&qSGL & 0.072 & 0.270 & 0.211 & 0.271 & 0.080 & 0.293 & 0.227 & 0.273 \\
		&4   &qL & 0.137 & 0.293 & 0.751 & 0.543 & 0.171 & 0.308 & 0.788 & 0.549 \\
		& 	&qGL & 0.333 & 0.618 & 0.306 & 0.397 & 0.337 & 0.630 & 0.319 & 0.435 \\
		\hline

		& 	&qSGL & 0.038 & 0.119 & 0.145 & 0.188 & 0.050 & 0.164 & 0.156 & 0.214 \\
		&1 &qL & 0.052 & 0.109 & 0.440 & 0.349 & 0.056 & 0.119 & 0.415 & 0.334 \\
		& 	&qGL &  0.307 & 0.501 & 0.229 & 0.333 & 0.316 & 0.562 & 0.279 & 0.400 \\

		\cline{2-11}

		& 	&qSGL & 0.036 & 0.100 & 0.141 & 0.173 & 0.046 & 0.146 & 0.157 & 0.202 \\
		&2&qL & 0.044 & 0.094 & 0.412 & 0.327 & 0.050 & 0.099 & 0.385 & 0.309 \\
		& 	&qGL& 0.305 & 0.498 & 0.227 & 0.330 & 0.316 & 0.560 & 0.278 & 0.413 \\
		\cline{2-11}

		400& &qSGL & 0.005 & 0.007 & 0.043 & 0.040 & 0.008 & 0.017 & 0.069 & 0.072 \\
		&3&qL & 0.007 & 0.007 & 0.076 & 0.059 & 0.009 & 0.013 & 0.154 & 0.121 \\
		&  &qGL & 0.291 & 0.430 & 0.198 & 0.301 & 0.294 & 0.445 & 0.209 & 0.312 \\
		\cline{2-11}

		& 	&Q & 0.028 & 0.080 & 0.135 & 0.160 & 0.038 & 0.122 & 0.150 & 0.191 \\
		&4&qL & 0.039 & 0.076 & 0.397 & 0.306 & 0.043 & 0.085 & 0.380 & 0.294 \\
		& 	&qGL & 0.302 & 0.485 & 0.223 & 0.325 & 0.311 & 0.532 & 0.263 & 0.388 \\
		\hline
	\end{tabular}
	\caption{Individual funtional $L_2$ error when SNR=10, as for Table \ref{tab:function error:snr5}.}
	%{ The first column n is the size of training data. The second column is the type of Noise. The third column is the method we used. ISE1: $||\hat{\beta}_1-\beta_1||_2^2$; ISE2: $||\hat{\beta}_2-\beta_2||_2^2$; ISE3: $||\hat{\beta}_3-\beta_3||_2^2$; ISE4: $||\hat{\beta}_4-\beta_4||_2^2$ }
	\label{tab:function error:snr10}
\end{table}

\begin{figure}[hbp]
	\centering
	\includegraphics[width=1\textwidth]{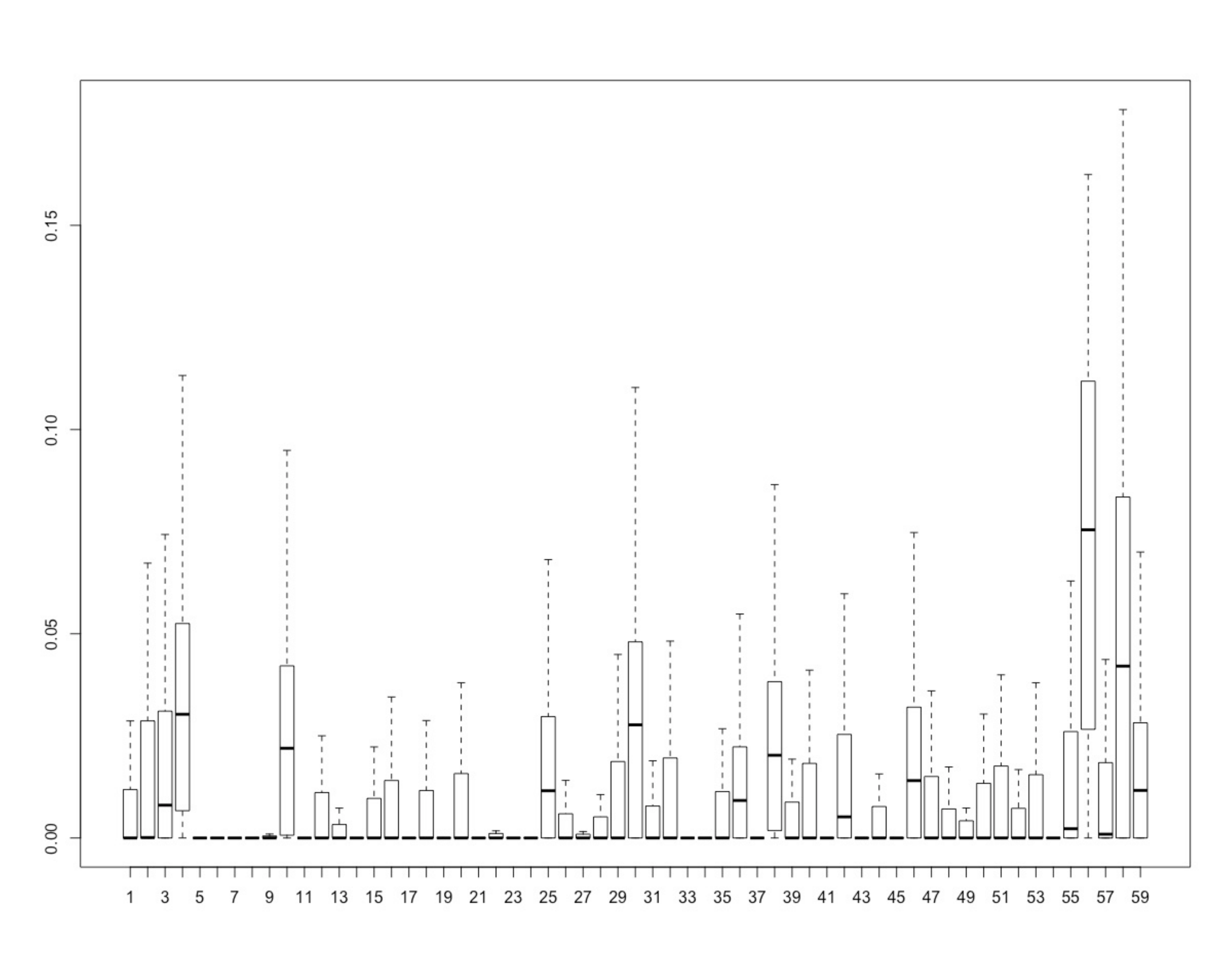}
	\caption{Boxplot of $L_2$ norm for each slope function, by using the quantile spare group lasso method.}
	\label{fig1:boxplot:sgl}
\end{figure}

\begin{figure}
	\centering
	\includegraphics[width=1\textwidth]{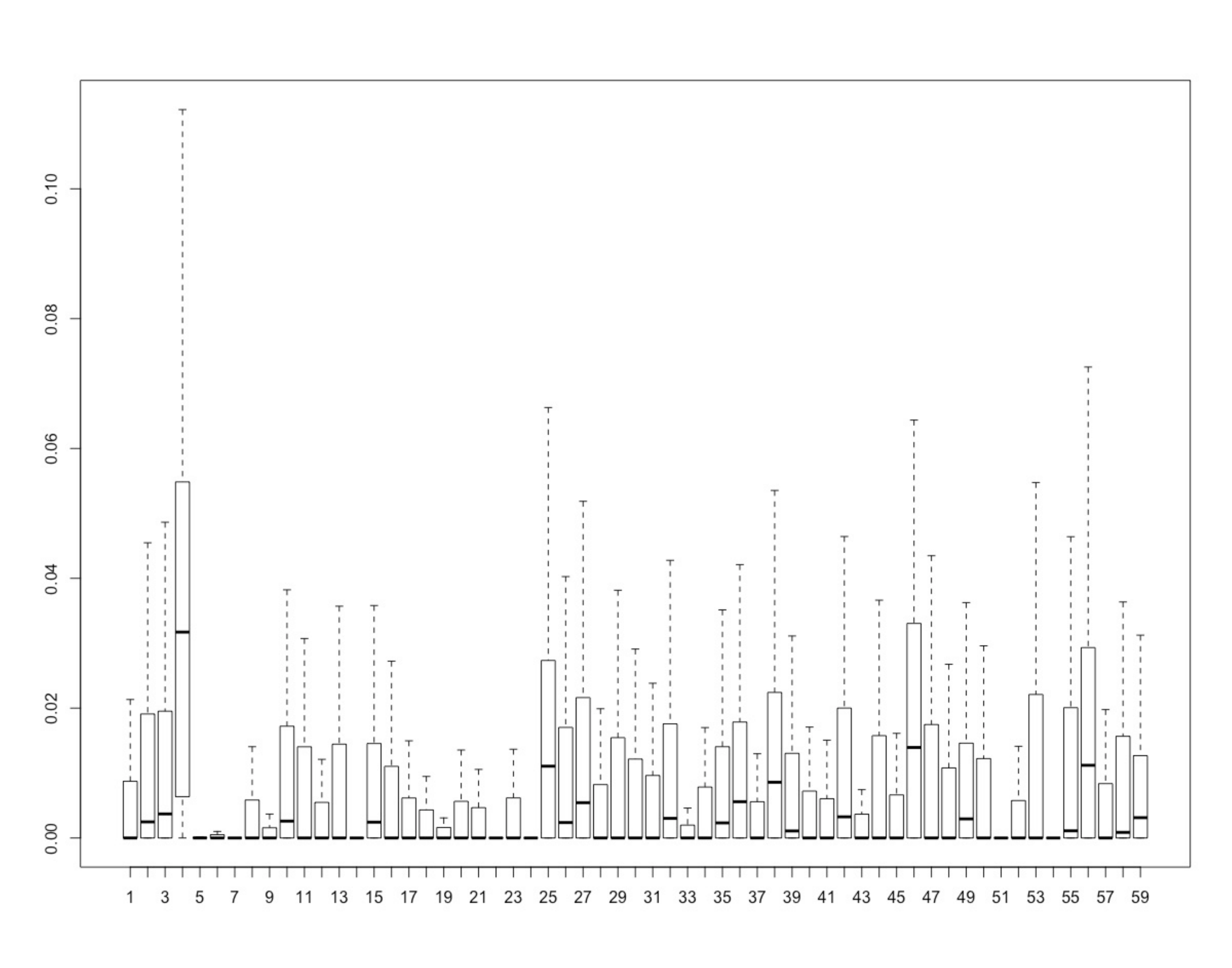}
	\caption{Boxplot of $L_2$ norm for each slope function, by using the quantile lasso method.}
	\label{fig1:boxplot:l}
\end{figure}

\begin{figure}
	\centering
	\includegraphics[width=1\textwidth]{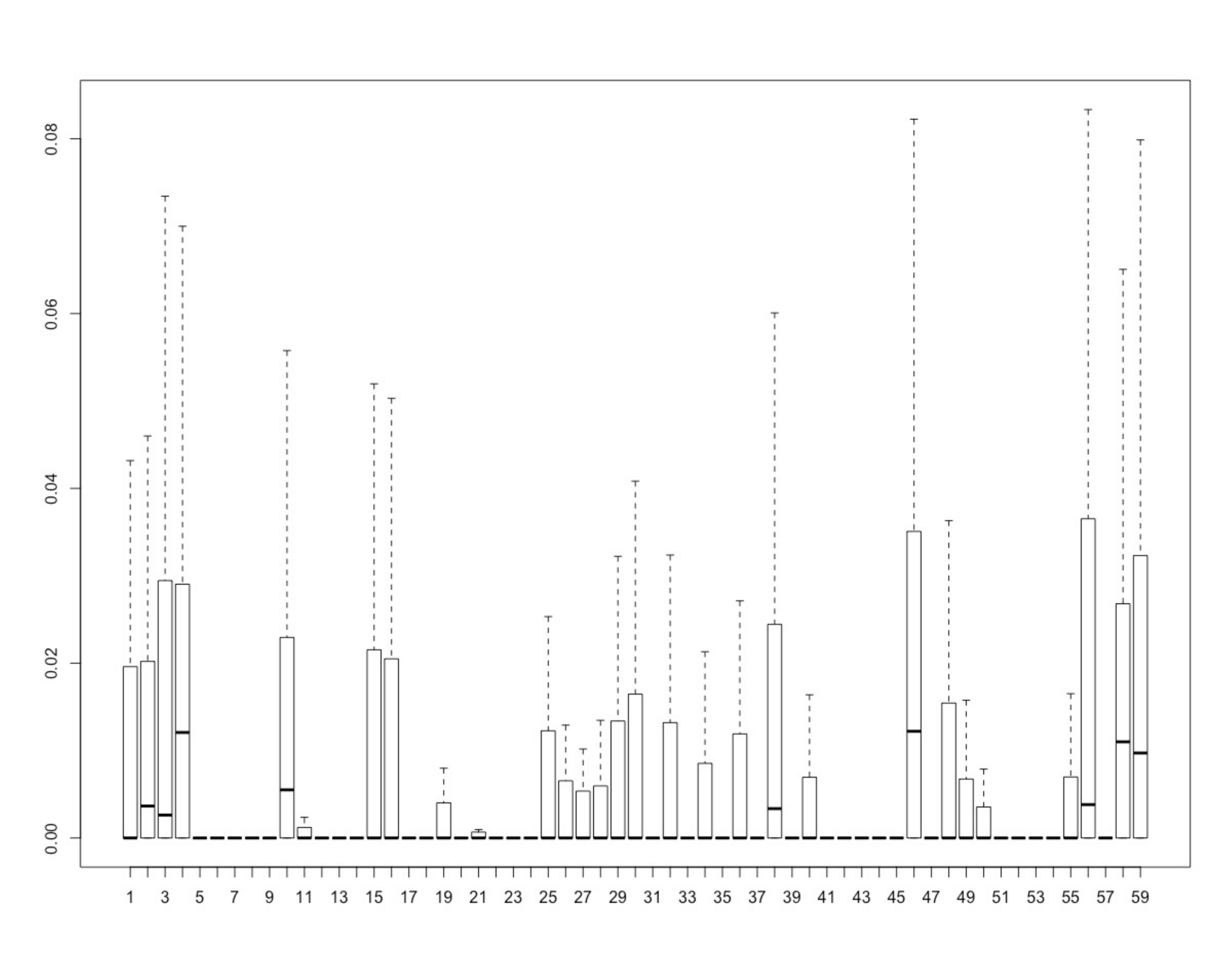}
	\caption{Boxplot of $L_2$ norm for each slope function, by using the quantile group lasso method.}
	\label{fig1:boxplot:g}
\end{figure}
\newpage
\subsection{Proof of Theorem 1}
%--------------------------------------  Proof, don?t delete, move the appendix ---------------------------------------------%
\begin{proof}
	First, we introduce some notation. The orthonormal wavelet basis set of $L^2[0,1]$ is defined as $\{\varphi_{j_0k},k=1,\dots,2^{j_0}\}\cup \{\psi_{jk},\ j\geq j_0,\ k=1,\dots,2^j\}$. Without loss of generality, the wavelet basis are ordered according to the scales from the coarsest level $J_0$ to the finest one. Let $\mathbb{V}_{N_n}:=\textrm{Span}\{\bm\varphi_1,\ \dots,\ \bm\varphi_{N_n}   \}$ be the space spanned by the first $N_n$ basis function, for example, if $N_n=2^{j_0+t}$, then the collection of  $\{\varphi_{j_0k},k=1\dots,2^{j_0}\}\cup\{\psi_{jk}, j_0\leq j\leq j_0+t-1,\ k=1,\dots,\ 2^j\}$ is the basis of $\mathbb{V}_{N_n}$. Let $\bm b_{N_n}^j$ be an $N_n\times1$ parameter vector with elements $b_k^j=\langle \beta_j(t),\varphi_{k}\rangle$. In addition, let $\beta_{N_n}^j$ be the functions reconstructed from the vector $\bm b_{N_n}^j$. Here $\beta_{N_n}^j$ is a linear approximation to $\beta_j$ by the first $N_n$ wavelet coefficients, while $\hat{\beta}_j$ denotes the function reconstructed from the wavelet coefficients $\hat{\bm b}_j$ from (\ref{cqr:sgl}).
	
	By the Parseval theorem, we have 	$||\hat{\beta}_j-\beta_j||_{L_2}^2=||\hat{b}^j_{N_n}-b^j_{N_n}||_2^2+\sum_{k={N_n}+1}^{\infty}{\theta_k^j}^2$. To derive the convergence rate of $\hat{\beta_j}$ to $\beta_j$, we bound the error in estimating $\beta_{N_n}^j$ by $\hat{\beta_j}$ and the error in approximating $\beta_j$ by $\beta_{N_n}$. By the Theorem 9.5 of \citet{Mallat:2008}, the linear approximation error goes to zero as
	\begin{equation}
		\sum_{k={N_n}+1}^{\infty}{b_k^j}^2=o(N_n^{-2d}).
	\end{equation}
	
	Let $\Upsilon^0=(\bm \alpha^0,\bm \gamma^0,\bm\theta^0)$ be the true coefficients with $\bm \theta^0=\rm{vec}^T \left(\bm b_{N_n}^1,\dots,\bm b_{N_n}^m \right)$.
	To obtain the result, we show that for any given $\varepsilon>0$, there exists a constant $C$ such that
	\begin{equation}\label{prop42}
		\textrm{Pr}\left\{\inf_{||\bm z||=C} L_n\left(\Upsilon^0+r_n \bm z\right)+ P_{\lambda_1,\lambda_2}\left(\bm\theta^0+r_n \bm z_{\theta}\right)>  L_n\left(\Upsilon^0\right)+ P_{\lambda_1,\lambda_2}\left(\bm\theta^0\right)\right\}\geq 1-\varepsilon,
	\end{equation}
	where $r_n=\sqrt{{N_n}/{n}}$ and $\bm z=(z_1,\dots,z_k,\bm z_{\bm\gamma},\bm z_{\bm\theta})$ is a vector with the same length of vector $\Upsilon^0$. This implies that there exists a local minimizer in the ball $\{\Upsilon^0+r_n z: ||z||\leq C\}$ with probability at least $1-\varepsilon$. Hence, there is a local minimizer $\widehat{\Upsilon}$ such that $||\widehat{\Upsilon}-\Upsilon^0||=O_p(r_n)$.
	%From this, we can also get $||\hat{\gamma}-\gamma^0||_2^2=O_p(r_n)$
	%		$|\hat{\alpha}_{\tau_k}-\alpha_{\tau_k}^0|=O_p(r_n)$ for $k=1,2 \dots s.$
	%	For any vector $\bm v=[z_1,\ \dots, z_k, \bm z_{\theta},\bm z_{u}]$ with $||\bm v||=C$, we need to compare the difference
	
	To show (\ref{prop42}), we compare $L_n(\bm\Upsilon^0)+P_n(\bm\theta^0)$ with $L_n(\bm\Upsilon^0+r_n \bm z)+P_n(\theta^0+r_n \bm z_{\theta})$.
	%	\begin{equation*}
	%L_n(\Upsilon^0+r_n v)-L_n(\Upsilon^0)+P_n(\bm\theta^0+r_n z_{\theta})-P_n(\theta^0).
	%	\end{equation*}
	By using the Knight identity,
	\begin{equation*}
		\rho_\tau(u-v)-\rho_\tau(u)=-v\varrho_\tau(u)+\int_{0}^{v}(I(u\leq t)-I(u\leq 0))dt,
	\end{equation*}
	where $\varrho_\tau(u)=\tau-I(u<0)$, we have
	\begin{eqnarray*}
		I&:=& L_n(\bm\Upsilon^0+r_n \bm v)-L_n(\bm\Upsilon^0)\\
		&=&\sum_{k=1}^{K} \sum_{i=1}^{n} [{\rho}_{\tau_k}(e_{ki}-d_{ki})- {\rho}_{\tau_k}(e_{ki})]\\
		&=&-\sum_{k=1}^{K} \sum_{i=1}^{n}[-d_{ki}\varrho_{\tau_k}(e_{ki})]+\sum_{k=1}^{K} \sum_{i=1}^{n}\int_{0}^{d_{ki}}(I(e_{ki}\leq t)-I(e_{ki}\leq 0))dt\\
		&=&I_1+I_2,
	\end{eqnarray*}
	where $e_{ki}=y_i-\alpha_{\tau_k}^0-\bm u_i^T\bm \gamma^0-\bm{v}_i^T\bm\theta^0$ and $b_{ki}=r_n z_k+r_n\bm u_i^T\bm z_u+r_n\bm{v}_i^T \bm z_\theta$. Note that $e_{ki}=\varepsilon_i-F^{-1}(\tau_k)+o(N_n^{-2d})$, hence we have $E(\varrho_{\tau_k}(e_{ki})=o(N_n^{-2d})$. By the definition of $d_{ki}$, we obtain $I_1 \leq {r_n}||\bm z||
	(\sum_{k=1}^{s} ||\sum_{i=1}^{n}\varrho_{\tau_k}(e_{ki})\bm A_i^T||)$
	and
	\begin{eqnarray*}
		E||\sum_{i=1}^{n}\varrho_{\tau_k}(e_{ki})\bm A_i||^2&=&E||\sum_{j=1}^{mN_n+1}\sum_{i=1}^{n}\sum_{l=1}^{n}a_{ij}a_{lj}\psi_{\tau_k}(e_{ki})\psi_{\tau_k}(e_{kl})||\\
		&=&O_p(nN_n),
	\end{eqnarray*}
	which leads to $E(I_1)\leq O_p( {r_n}\sqrt{nN_n})||\bm z||=O_p(nr_n^2)||\bm z||$.
	
	Now, we consider the expectation of $I_2$. Using the expression of $e_{ki}$, we get
	\begin{eqnarray*}
		E(I_2)& =& \sum_{k=1}^{K} \sum_{i=1}^{n}\int_{0}^{d_{ki}}(\textrm{Pr}(e_{ki}\leq t)-\textrm{Pr}(e_{ki}\leq 0))dt\\
		&=&\sum_{k=1}^{K} \sum_{i=1}^{n}\int_{0}^{d_{ki}}(F(F^{-1}(\tau_k)+o(N_n^{-2d})+t)-F(F^{-1}(\tau_k)+o(N_n^{-2d})))dt\\
		&=&\sum_{k=1}^{K} \sum_{i=1}^{n}\int_{0}^{d_{ki}}(f(F^{-1}(\tau_k)+o(N_n^{-2d}))t+\frac {f'(\xi)} {2} t^2)dt,
	\end{eqnarray*}
	where $\xi$ lies between $F^{-1}(\tau_k)+o(N_n^{-2d})$ and $F^{-1}(\tau_k)+o(N_n^{-2d})+d_{ki}$. Since there exists $M$ such that $||A_i||_2^2<M$, we have
	\begin{equation*}
		\max_{1 \leq i\leq n} |r_nz_k+r_n\bm v_i^T\bm z_\theta|\to 0.
	\end{equation*}
	Then,  the lower bound of $E(I_2)$ is of the form
	\begin{eqnarray*}
		E(I_2)&=&\frac 1 {2} r_n^2\sum_{k=1}^{K}\{[f(F^{-1}(\tau_k)+o(N_n^{-2d}))+o_p(1)](\bm g_k^TA^TA\bm g_k)\}\\
		&\geq& \frac{c_1n r_n^2}{2}||\bm z||^2_2\min_k\{f(F^{-1}(\tau_k)+o(N_n^{-2d}))+o_p(1) \},
	\end{eqnarray*}
	where $\bm g_k$ is a vector, such as $\bm g_k=(z_k,\bm z_{\theta}^T,\bm z_{u}^T)^T$. Finally, since $r_n \to 0$ and $||\bm z||_2\leq C$, we have
	\begin{eqnarray*}
		II:= P_n(\theta^0+r_n \bm z_\theta)-P_n(\theta^0)&\leq& \lambda_1r_n||\bm z_\theta||_{1}+\lambda_2 r_n \sum_{j=1}^{m}||\bm z_{\theta_j}||_2\\
		&\leq& \lambda_1 r_n \sqrt{mN}||\bm z_\theta||_2+\lambda_2r_nm||\bm z_\theta||_2\\
		&= & O_p(nr_n^2||\bm z_\theta||_2).
	\end{eqnarray*}
	Since $II$ is bounded by $r_n^2||\bm z_\theta||_2$ ,we can choose a $C$ such that the II is dominated by the term $I_2$ on $||u||=C$ uniformly. So  $Q_n(\Sigma^0+r_n u)- Q_n(\Sigma^0)>0$ holds uniformly on $||u||=C$. This completes the proof.
\end{proof}
%--------------------------------------  Proof move the appendix ---------------------------------------------%

\bibliographystyle{chicago}
\bibliography{refsg}
\end{document}